\documentclass{amsart}
\usepackage{amssymb}
\newcommand{\BC}{\mathbb C}

\newcommand{\BZ}{\mathbb Z}

\newcommand{\cO}{{\mathcal O}}
\newcommand{\cL}{{\mathcal L}}
\newcommand{\cN}{{\mathcal N}}
\newcommand{\uc}{\bf c}
\newcommand{\bc}{\uc}

\newcommand{\fU}{{\mathfrak U}}
\newcommand{\fg}{{\mathfrak g}}
\newcommand{\fh}{{\mathfrak h}}

\newcommand{\sq}{$\square$}
\def\bbb#1#2#3{\mbox{\ \begin{picture}(20,5)(-5,0)
    \put(0,0){\circle*{2}}
    \put(- 2,-8){\mbox{$\scriptstyle #1$}}
    \put(0,0){\line(1,0){10}}
    \put(10,0){\circle*{2}}
    \put(8,-8){\mbox{$\scriptstyle #2$}}
    \put(13,-2){\mbox{$\scriptstyle >$}}
    \bezier{30}(10,0)(15,2)(20,0)
    \bezier{30}(10,0)(15,-2)(20,0)
    \put(20,0){\circle*{2}}
    \put(18,-8){\mbox{$\scriptstyle #3$}}
  \end{picture}\ \makebox(0,15){}}}
\def\bbbb#1#2#3#4{\mbox{\ \begin{picture}(30,5)(-5,0)
    \put(0,0){\circle*{2}}
    \put(- 2,-8){\mbox{$\scriptstyle #1$}}
    \put(0,0){\line(1,0){10}}
    \put(10,0){\circle*{2}}
    \put(8,-8){\mbox{$\scriptstyle #2$}}
    \put(10,0){\line(1,0){10}}
    \put(20,0){\circle*{2}}
    \put(18,-8){\mbox{$\scriptstyle #3$}}
    \put(23,-2){\mbox{$\scriptstyle >$}}
    \bezier{30}(20,0)(25,2)(30,0)
    \bezier{30}(20,0)(25,-2)(30,0)
    \put(30,0){\circle*{2}}
    \put(28,-8){\mbox{$\scriptstyle #4$}}
  \end{picture}\ \makebox(0,15){}}}
\def\bb#1#2{\mbox{\ \begin{picture}(10,5)(-5,0)
    \put(0,0){\circle*{2}}
    \put(- 2,-8){\mbox{$\scriptstyle #1$}}
    \put(2,-2){\mbox{$>$}}
    \bezier{30}(0,0)(5,2)(10,0)
    \bezier{30}(0,0)(5,-2)(10,0)
    \put(10,0){\circle*{2}}
    \put(8,-8){\mbox{$\scriptstyle #2$}}
  \end{picture}\ \makebox(0,15){}}}
\def\bbbbb#1#2#3#4#5{\mbox{\ \begin{picture}(40,5)(-5,0)
    \put(0,0){\circle*{2}}
    \put(- 2,-8){\mbox{$\scriptstyle #1$}}
    \put(0,0){\line(1,0){10}}
    \put(10,0){\circle*{2}}
    \put(8,-8){\mbox{$\scriptstyle #2$}}
    \put(10,0){\line(1,0){10}}
    \put(20,0){\circle*{2}}
    \put(18,-8){\mbox{$\scriptstyle #3$}}
    \put(20,0){\line(1,0){10}}
    \put(30,0){\circle*{2}}
    \put(28,-8){\mbox{$\scriptstyle #4$}}
    \put(33,-2){\mbox{$\scriptstyle >$}}
    \bezier{30}(30,0)(35,2)(40,0)
    \bezier{30}(30,0)(35,-2)(40,0)
    \put(40,0){\circle*{2}}
    \put(38,-8){\mbox{$\scriptstyle #5$}}
  \end{picture}\ \makebox(0,15){}}}
\def\bbbbbb#1#2#3#4#5#6{\mbox{\ \begin{picture}(50,5)(-5,0)
    \put(0,0){\circle*{2}}
    \put(- 2,-8){\mbox{$\scriptstyle #1$}}
    \put(0,0){\line(1,0){10}}
    \put(10,0){\circle*{2}}
    \put(8,-8){\mbox{$\scriptstyle #2$}}
    \put(10,0){\line(1,0){10}}
    \put(20,0){\circle*{2}}
    \put(18,-8){\mbox{$\scriptstyle #3$}}
    \put(20,0){\line(1,0){10}}
    \put(30,0){\circle*{2}}
    \put(28,-8){\mbox{$\scriptstyle #4$}}
    \put(30,0){\line(1,0){10}}
    \put(40,0){\circle*{2}}
    \put(38,-8){\mbox{$\scriptstyle #5$}}
    \put(43,-2){\mbox{$\scriptstyle >$}}
    \bezier{30}(40,0)(45,2)(50,0)
    \bezier{30}(40,0)(45,-2)(50,0)
    \put(50,0){\circle*{2}}
    \put(48,-8){\mbox{$\scriptstyle #6$}}
  \end{picture}\ \makebox(0,15){}}}
\def\bbbbbbb#1#2#3#4#5#6#7{\mbox{\ \begin{picture}(60,5)(-5,0)
    \put(0,0){\circle*{2}}
    \put(- 2,-8){\mbox{$\scriptstyle #1$}}
    \put(0,0){\line(1,0){10}}
    \put(10,0){\circle*{2}}
    \put(8,-8){\mbox{$\scriptstyle #2$}}
    \put(10,0){\line(1,0){10}}
    \put(20,0){\circle*{2}}
    \put(18,-8){\mbox{$\scriptstyle #3$}}
    \put(20,0){\line(1,0){10}}
    \put(30,0){\circle*{2}}
    \put(28,-8){\mbox{$\scriptstyle #4$}}
    \put(30,0){\line(1,0){20}}
    \put(40,0){\circle*{2}}
    \put(38,-8){\mbox{$\scriptstyle #5$}}
    \put(40,0){\line(1,0){10}}
    \put(50,0){\circle*{2}}
    \put(48,-8){\mbox{$\scriptstyle #6$}}
    \put(52,-2){\mbox{$>$}}
    \bezier{30}(50,0)(55,2)(60,0)
    \bezier{30}(50,0)(55,-2)(60,0)
    \put(60,0){\circle*{2}}
    \put(58,-8){\mbox{$\scriptstyle #7$}}
  \end{picture}\ \makebox(0,15){}}}
\def\sp#1#2#3{\mbox{\ \begin{picture}(20,5)(0,0)
    \put(0,0){\circle*{2}}
    \put(- 2,-8){\mbox{$\scriptstyle #1$}}
    \put(0,0){\line(1,0){10}}
    \put(10,0){\circle*{2}}
    \put(8,-8){\mbox{$\scriptstyle #2$}}
    \put(12,-2){\mbox{$\scriptstyle <$}}
    \bezier{30}(10,0)(15,2)(20,0)
    \bezier{30}(10,0)(15,-2)(20,0)
    \put(20,0){\circle*{2}}
    \put(18,-8){\mbox{$\scriptstyle #3$}}
  \end{picture}\ \makebox(0,15){}}}
\def\spp#1#2#3#4{\mbox{\ \begin{picture}(30,5)(0,0)
    \put(0,0){\circle*{2}}
    \put(- 2,-8){\mbox{$\scriptstyle #1$}}
    \put(0,0){\line(1,0){10}}
    \put(10,0){\circle*{2}}
    \put(8,-8){\mbox{$\scriptstyle #2$}}
    \put(10,0){\line(1,0){10}}
    \put(20,0){\circle*{2}}
    \put(18,-8){\mbox{$\scriptstyle #3$}}
    \put(22,-2){\mbox{$\scriptstyle <$}}
    \bezier{30}(20,0)(25,2)(30,0)
    \bezier{30}(20,0)(25,-2)(30,0)
    \put(30,0){\circle*{2}}
    \put(28,-8){\mbox{$\scriptstyle #4$}}
  \end{picture}\ \makebox(0,15){}}}
    
\def\s#1#2{\mbox{\ \begin{picture}(10,5)(0,0)
    \put(0,0){\circle*{2}}
    \put(- 2,-8){\mbox{$\scriptstyle #1$}}
    \put(2,-2){\mbox{$<$}}
    \bezier{30}(0,0)(5,2)(10,0)
    \bezier{30}(0,0)(5,-2)(10,0)
    \put(10,0){\circle*{2}}
    \put(8,-8){\mbox{$\scriptstyle #2$}}
  \end{picture}\ \makebox(0,15){}}}
\def\sppp#1#2#3#4#5{\mbox{\ \begin{picture}(40,5)(0,0)
    \put(0,0){\circle*{2}}
    \put(- 2,-8){\mbox{$\scriptstyle #1$}}
    \put(0,0){\line(1,0){10}}
    \put(10,0){\circle*{2}}
    \put(8,-8){\mbox{$\scriptstyle #2$}}
    \put(10,0){\line(1,0){10}}
    \put(20,0){\circle*{2}}
    \put(18,-8){\mbox{$\scriptstyle #3$}}
    \put(20,0){\line(1,0){10}}
    \put(30,0){\circle*{2}}
    \put(28,-8){\mbox{$\scriptstyle #4$}}
    \put(32,-2){\mbox{$\scriptstyle <$}}
    \bezier{30}(30,0)(35,2)(40,0)
    \bezier{30}(30,0)(35,-3)(40,0)
    \put(40,0){\circle*{2}}
    \put(38,-8){\mbox{$\scriptstyle #5$}}
  \end{picture}\ \makebox(0,15){}}}
\def\ssp#1#2#3#4#5#6{\mbox{\ \begin{picture}(50,5)(0,0)
    \put(0,0){\circle*{2}}
    \put(- 2,-8){\mbox{$\scriptstyle #1$}}
    \put(0,0){\line(1,0){10}}
    \put(10,0){\circle*{2}}
    \put(8,-8){\mbox{$\scriptstyle #2$}}
    \put(10,0){\line(1,0){10}}
    \put(20,0){\circle*{2}}
    \put(18,-8){\mbox{$\scriptstyle #3$}}
    \put(20,0){\line(1,0){10}}
    \put(30,0){\circle*{2}}
    \put(28,-8){\mbox{$\scriptstyle #4$}}
    \put(30,0){\line(1,0){10}}
    \put(40,0){\circle*{2}}
    \put(38,-8){\mbox{$\scriptstyle #5$}}
    \put(42,-2){\mbox{$\scriptstyle <$}}
    \bezier{30}(40,0)(45,2)(50,0)
    \bezier{30}(40,0)(45,-2)(50,0)
    \put(50,0){\circle*{2}}
    \put(48,-8){\mbox{$\scriptstyle #6$}}
  \end{picture}\ \makebox(0,15){}}}
\def\sssp#1#2#3#4#5#6#7{\mbox{\ \begin{picture}(60,5)(0,0)
    \put(0,0){\circle*{2}}
    \put(- 2,-8){\mbox{$\scriptstyle #1$}}
    \put(0,0){\line(1,0){10}}
    \put(10,0){\circle*{2}}
    \put(8,-8){\mbox{$\scriptstyle #2$}}
    \put(10,0){\line(1,0){10}}
    \put(20,0){\circle*{2}}
    \put(18,-8){\mbox{$\scriptstyle #3$}}
    \put(20,0){\line(1,0){10}}
    \put(30,0){\circle*{2}}
    \put(28,-8){\mbox{$\scriptstyle #4$}}
    \put(30,0){\line(1,0){10}}
    \put(40,0){\circle*{2}}
    \put(38,-8){\mbox{$\scriptstyle #5$}}
    \put(40,0){\line(1,0){10}}\
    \put(50,0){\circle*{2}}
    \put(48,-8){\mbox{$\scriptstyle #6$}}
    \put(51,-2){\mbox{$<$}}
    \bezier{30}(50,0)(55,2)(60,0)
    \bezier{30}(50,0)(55,-2)(60,0)
    \put(60,0){\circle*{2}}
    \put(58,-8){\mbox{$\scriptstyle #7$}}
  \end{picture}\ \makebox(0,15){}}}
\def\so#1#2#3{\mbox{\ \begin{picture}(10,6)(0,0)
    \put(0,3){\circle*{2}}
    \put(- 2,-5){\mbox{$\scriptstyle #1$}}
    \put(0,3){\line(5,1){10}}
    \put(10,5){\circle*{2}}
    \put(0,3){\line(5,-1){10}}
    \put(10,1){\circle*{2}}
    \put(12,3){\mbox{$\scriptstyle #2$}}
    \put(12,-3){\mbox{$\scriptstyle #3$}}
    \end{picture}\ \makebox(0,20){}}}
\def\soo#1#2#3#4{\mbox{\ \begin{picture}(20,6)(0,0)
    \put(0,3){\circle*{2}}
    \put(- 2,-5){\mbox{$\scriptstyle #1$}}
    \put(0,3){\line(1,0){10}}
    \put(10,3){\circle*{2}}
    \put(8,-5){\mbox{$\scriptstyle #2$}}
    \put(10,3){\line(5,1){10}}
    \put(20,5){\circle*{2}}
    \put(10,3){\line(5,-1){10}}
    \put(20,1){\circle*{2}}
    \put(22,3){\mbox{$\scriptstyle #3$}}
    \put(22,-3){\mbox{$\scriptstyle #4$}}
    \end{picture}\ \makebox(0,20){}}}
\def\sooo#1#2#3#4#5{\mbox{\ \begin{picture}(30,6)(0,0)
    \put(0,3){\circle*{2}}
    \put(- 2,-5){\mbox{$\scriptstyle #1$}}
    \put(0,3){\line(1,0){10}}
    \put(10,3){\circle*{2}}
    \put(8,-5){\mbox{$\scriptstyle #2$}}
    \put(10,3){\line(1,0){10}}
    \put(20,3){\circle*{2}}
    \put(18,-5){\mbox{$\scriptstyle #3$}}
    \put(20,3){\line(5,1){10}}
    \put(30,5){\circle*{2}}
    \put(20,3){\line(5,-1){10}}
    \put(30,1){\circle*{2}}
    \put(32,3){\mbox{$\scriptstyle #4$}}
    \put(32,-3){\mbox{$\scriptstyle #5$}}
    \end{picture}\ \makebox(0,20){}}}
\def\sso#1#2#3#4#5#6{\mbox{\ \begin{picture}(45,6)(0,0)
    \put(0,3){\circle*{2}}
    \put(- 2,-5){\mbox{$\scriptstyle #1$}}
    \put(0,3){\line(1,0){10}}
    \put(10,3){\circle*{2}}
    \put(8,-5){\mbox{$\scriptstyle #2$}}
    \put(10,3){\line(1,0){10}}
    \put(20,3){\circle*{2}}
    \put(18,-5){\mbox{$\scriptstyle #3$}}
    \put(20,3){\line(1,0){10}}
    \put(30,3){\circle*{2}}
    \put(28,-5){\mbox{$\scriptstyle #4$}}
    \put(30,3){\line(5,1){10}}
    \put(40,5){\circle*{2}}
    \put(30,3){\line(5,-1){10}}
    \put(40,1){\circle*{2}}
    \put(42,3){\mbox{$\scriptstyle #5$}}
    \put(42,-3){\mbox{$\scriptstyle #6$}}
    \end{picture}\ \makebox(0,20){}}}
\def\ssso#1#2#3#4#5#6#7{\mbox{\ \begin{picture}(55,6)(0,0)
    \put(0,3){\circle*{2}}
    \put(- 2,-5){\mbox{$\scriptstyle #1$}}
    \put(0,3){\line(1,0){10}}
    \put(10,3){\circle*{2}}
    \put(8,-5){\mbox{$\scriptstyle #2$}}
    \put(10,3){\line(1,0){10}}
    \put(20,3){\circle*{2}}
    \put(18,-5){\mbox{$\scriptstyle #3$}}
    \put(20,3){\line(1,0){10}}
    \put(30,3){\circle*{2}}
    \put(28,-5){\mbox{$\scriptstyle #4$}}
    \put(30,3){\line(1,0){10}}
    \put(40,3){\circle*{2}}
    \put(38,-5){\mbox{$\scriptstyle #5$}}
    \put(40,3){\line(5,1){10}}
    \put(50,5){\circle*{2}}
    \put(40,3){\line(5,-1){10}}
    \put(50,1){\circle*{2}}
    \put(52,3){\mbox{$\scriptstyle #6$}}
    \put(52,-3){\mbox{$\scriptstyle #7$}}
    \end{picture}\ \makebox(0,20){}}}
\def\ssoo#1#2#3#4#5#6#7#8{\mbox{\
\begin{picture}(60,6)(0,0)
    \put(0,3){\circle*{2}}
    \put(- 2,-5){\mbox{$\scriptstyle #1$}}
    \put(0,3){\line(1,0){10}}
    \put(10,3){\circle*{2}}
    \put(8,-5){\mbox{$\scriptstyle #2$}}
    \put(10,3){\line(1,0){10}}
    \put(20,3){\circle*{2}}
    \put(18,-5){\mbox{$\scriptstyle #3$}}
    \put(20,3){\line(1,0){10}}
    \put(30,3){\circle*{2}}
    \put(28,-5){\mbox{$\scriptstyle #4$}}
    \put(30,3){\line(1,0){10}}
    \put(40,3){\circle*{2}}
    \put(38,-5){\mbox{$\scriptstyle #5$}}
    \put(40,3){\line(1,0){10}}
    \put(50,3){\circle*{2}}
    \put(48,-5){\mbox{$\scriptstyle #6$}}
    \put(50,3){\line(5,1){10}}
    \put(60,5){\circle*{2}}
    \put(50,3){\line(5,-1){10}}
    \put(60,1){\circle*{2}}
    \put(62,3){\mbox{$\scriptstyle #7$}}
    \put(62,-3){\mbox{$\scriptstyle #8$}}
    \end{picture}\ \makebox(0,20){}}}
\def\g#1#2{\mbox{\ \begin{picture}(20,5)(0,0)
    \put(0,0){\circle*{2}}
    \put(- 2,-8){\mbox{$\scriptstyle #1$}}
    \put(6,-2){\mbox{$<$}}
    \bezier{60}(0,0)(10,4)(20,0)
    \put(0,0){\line(1,0){20}}
    \bezier{60}(0,0)(10,-4)(20,0)
    \put(20,0){\circle*{2}}
    \put(18,-8){\mbox{$\scriptstyle #2$}}
  \end{picture}\ \makebox(0,15){}}}
\def\f#1#2#3#4{\mbox{\ \begin{picture}(30,5)(0,0)
    \put(0,0){\circle*{2}}
    \put(- 2,-8){\mbox{$\scriptstyle #1$}}
    \put(0,0){\line(1,0){10}}
    \put(10,0){\circle*{2}}
    \put(8,-8){\mbox{$\scriptstyle #2$}}
    \put(13,-2){\mbox{$>$}}
    \bezier{30}(10,0)(15,2)(20,0)
    \bezier{30}(10,0)(15,-2)(20,0)
    \put(20,0){\circle*{2}}
    \put(18,-8){\mbox{$\scriptstyle #3$}}
    \put(20,0){\line(1,0){10}}
    \put(30,0){\circle*{2}}
    \put(28,-8){\mbox{$\scriptstyle #4$}}
\end{picture}\ \makebox(0,15){}}}
\def\e#1#2#3#4#5#6{\mbox{\ \begin{picture}(40,10)(0,0)
    \put(0,0){\circle*{2}}
    \put(- 2,-8){\mbox{$\scriptstyle #1$}}
    \put(0,0){\line(1,0){10}}
    \put(10,0){\circle*{2}}
    \put(8,-8){\mbox{$\scriptstyle #3$}}
    \put(10,0){\line(1,0){10}}
    \put(20,0){\circle*{2}}
    \put(18,-8){\mbox{$\scriptstyle #4$}}
    \put(20,0){\line(1,0){10}}
    \put(30,0){\circle*{2}}
    \put(28,-8){\mbox{$\scriptstyle #5$}}
    \put(30,0){\line(1,0){10}}
    \put(40,0){\circle*{2}}
    \put(38,-8){\mbox{$\scriptstyle #6$}}
    \put(20,10){\circle*{2}}
    \put(20,10){\line(0,-1){10}}
    \put(22,7){\mbox{$\scriptstyle #2$}}
\end{picture}\ \makebox(0,35){}}}
\def\ee#1#2#3#4#5#6#7{\mbox{\ \begin{picture}(50,10)(0,0)
    \put(0,0){\circle*{2}}
    \put(- 2,-8){\mbox{$\scriptstyle #1$}}
    \put(0,0){\line(1,0){10}}
    \put(10,0){\circle*{2}}
    \put(8,-8){\mbox{$\scriptstyle #3$}}
    \put(10,0){\line(1,0){10}}
    \put(20,0){\circle*{2}}
    \put(18,-8){\mbox{$\scriptstyle #4$}}
    \put(20,0){\line(1,0){10}}
    \put(30,0){\circle*{2}}
    \put(28,-8){\mbox{$\scriptstyle #5$}}
    \put(30,0){\line(1,0){10}}
    \put(40,0){\circle*{2}}
    \put(38,-8){\mbox{$\scriptstyle #6$}}
    \put(40,0){\line(1,0){10}}
    \put(50,0){\circle*{2}}
    \put(48,-8){\mbox{$\scriptstyle #7$}}
    \put(20,10){\circle*{2}}
    \put(20,10){\line(0,-1){10}}
    \put(22,7){\mbox{$\scriptstyle #2$}}
\end{picture}\ \makebox(0,35){}}}
\def\eee#1#2#3#4#5#6#7#8{\mbox{\ \begin{picture}(60,10)(0,0)
    \put(0,0){\circle*{2}}
    \put(- 2,-8){\mbox{$\scriptstyle #1$}}
    \put(0,0){\line(1,0){10}}
    \put(10,0){\circle*{2}}
    \put(8,-8){\mbox{$\scriptstyle #3$}}
    \put(10,0){\line(1,0){10}}
    \put(20,0){\circle*{2}}
    \put(18,-8){\mbox{$\scriptstyle #4$}}
    \put(20,0){\line(1,0){10}}
    \put(30,0){\circle*{2}}
    \put(28,-8){\mbox{$\scriptstyle #5$}}
    \put(30,0){\line(1,0){10}}
    \put(40,0){\circle*{2}}
    \put(38,-8){\mbox{$\scriptstyle #6$}}
    \put(40,0){\line(1,0){10}}
    \put(50,0){\circle*{2}}
    \put(48,-8){\mbox{$\scriptstyle #7$}}
    \put(50,0){\line(1,0){10}}
    \put(60,0){\circle*{2}}
    \put(58,-8){\mbox{$\scriptstyle #8$}}
    \put(20,10){\circle*{2}}
    \put(20,10){\line(0,-1){10}}
    \put(22,7){\mbox{$\scriptstyle #2$}}
\end{picture}\ \makebox(0,35){}}}

\title{Calculating canonical distinguished involutions in the affine Weyl
groups}
\author{Tanya Chmutova} \email{chmutova@mccme.ru} \address{Independent Moscow
University, Bolshoj Vlas'evskij per., 11, Moscow 121002 Russia}
\author{Viktor Ostrik} \email{ostrik@math.mit.edu} \address{Department of
Mathematics, M.I.T., Cambridge, MA 02139 USA}
\date{June 2001}
\thanks{The second author was partially supported by
an NSF grant}

\begin{document}
\begin{abstract}
Distinguished involutions in the affine Weyl groups, defined 
by G.~Lusztig, play an essential role in the Kazhdan-Lusztig combinatorics of 
these groups. A distinguished involution is called canonical if it is the 
shortest element in its double coset with respect to the finite Weyl group. 
Each two-sided cell in the affine Weyl group contains precisely one canonical
distinguished involution. In this note we calculate the canonical 
distinguished involutions in the affine Weyl groups of rank $\le 7$. We also 
prove some partial results relating canonical distinguished involutions and 
Dynkin's diagrams of the nilpotent orbits in the Langlands dual group.
\end{abstract} 
\maketitle

\section{Introduction}
It is well known that the Kazhdan-Lusztig combinatorics of the affine Hecke
algebra is deeply related with the geometry of the corresponding algebraic
group (over the complex numbers) $G$. The central result here is the deep
Theorem of George Lusztig
establishing a bijection between the set $\fU$ of unipotent classes in $G$ and
the set of {\em two-sided cells} in the corresponding affine Weyl group $W_a$,
see
\cite{L4}. Using this one defines a map from $\fU$ to the set $X_+$ of the
dominant weights in the following way: let $\cO$ be an unipotent orbit and
let $\bc_{\cO}$ be the corresponding two-sided cell. George Lusztig and
Nanhua Xi attached to $\bc_{\cO}$ a {\em canonical left cell} $C_{\cO}$, see
\cite{LX}. In turn the left cell $C_{\cO}$ contains a unique {\em distinguished
involution} $d_{\cO}\in W_a$, see \cite{L2}, which is the shortest element in
its double coset $Wd_\cO W$ with respect to the finite Weyl group $W\subset
W_a$. It is well known that the set of double cosets $W\setminus W_a/W$ is 
bijective to the set $X_+$ since any such coset contains unique translation by
dominant weight. Combining maps above we get a canonical map $\cL :\fU \to 
X_+$. The explicit calculation of this map is equivalent to the determination 
of the distinguished involutions lying in the canonical cells (or, 
equivalently, which are shortest in their left coset with respect to the 
finite Weyl group). We call such involutions {\em canonical distinguished 
involutions}. We believe that understanding of these involutions is an 
important step towards understanding of all distinguished involutions and 
cells in the affine Weyl group.

In \cite{O} one of us suggested a conjectural algorithm for the calculation
of $\cL$ and now this algorithm is known to be correct thanks to the
(unfortunately still undocumented) work of R.~Bezrukavnikov. The aim of this
note is to present results of calculations using this algorithm.

The paper is organized as follows. In the section 2 we recollect necessary
facts. In section 3 we present our main results: calculation of the map
$\cL$ for group $GL_n$ and partial results for other groups. These results
seems to be known to the experts but to the best of our knowledge were never
published. In section 4 we present tables with the results of explicit 
calculation of the map $\cL$ for groups of small rank. These tables should be 
considered as a main result of this work.

We would like to thank David Vogan and Pramod Achar for useful conversations.

\section{Recollections}
\subsection{Notations}
Let $G$ be a semisimple algebraic group over the complex numbers. Let $\fg$ 
denote the Lie algebra of $G$ and let $\cN \subset \fg$ denote
the nilpotent cone. As a $G-$variety $\cN$ is isomorphic to the subvariety of
unipotent elements of $G$ via the exponential map but it has one additional 
virtue --- an obvious action of $\BC^*$ by dilations commuting with the 
$G-$action. We will consider $\cN$ as $G\times \BC^*-$variety via the 
following action $(g,z)n=z^{-2}Ad(g)n$ for $(g,z)\in G\times \BC^*$ and 
$n\in \cN$.

The variety $\cN$ consists of finitely many $G-$orbits, see e.g. \cite{CM}.
These orbits called {\em nilpotent orbits} are the main subject of our study.
Any nilpotent orbit $\cO$ is identified via its {\em Dynkin diagram} defined as
follows: let $e\in \cO$ be a representative, by the Jacobson---Morozov Theorem
it can be included in $sl_2-$triple $(e, f, h)$ (i.e. $[h, e]=2e,\; 
[h, f]=-2f,\; [e, f]=h$). The semisimple element $h$ is uniquely defined up to
$G-$conjugacy by $\cO$ and vice versa, see e.g. \cite{CM}. Let $\fh \subset 
\fg$ be a Cartan subalgebra, let $R\subset \fh^*$ be the root system, choose a
subset $R_+\subset R$ of positive roots and let $\{ \alpha_i, i\in I\}$ be the
set of simple roots ($I$ is the set of vertices of Dynkin diagram of $\fg$).
The element $h$ is conjugate to a unique $h_0\in \fh$ such that
$\alpha_i(h_0)$ is positive for any $i\in I$ and moreover $\alpha_i(h_0)\in
\{ 0, 1, 2\}$, see \cite{CM}. Thus any nilpotent orbit can be specialized
by labeling the Dynkin diagram of $\fg$ by numbers $0, 1, 2$ and this is called
the (labeled) Dynkin diagram of $\fg$. We note that the Dynkin diagram $h_0$ is
naturally integral dominant {\em coweight} for group $G$.

Let $X$ denote the weight lattice of $G$ and let $X_+$ denote the set of
dominant weights.

\subsection{Review of \cite{O}}
Let $K_{G\times \BC^*}(\cN )$ denote the Grothendieck group of the category
of $G\times \BC^*-$equivariant coherent sheaves on $\cN$. It has an obvious
structure of a module over representation ring $Rep(\BC^*)$ of $\BC^*$. Let
$v$ denote tautological representation of $\BC^*$. Then $Rep(\BC^*)=\BZ
[v, v^{-1}]$ and the rule $v\mapsto v^{-1}$ defines an involution
$\bar : Rep(\BC^*)\to Rep(\BC^*)$.

In \cite{O} it was constructed the following:

1) the basis $\{ AJ(\lambda )\}$ of $K_{G\times \BC^*}(\cN )$ over $Rep(\BC^*)=
\BZ [v, v^{-1}]$ labeled by dominant weights $\lambda \in X_+$.

2) $Rep(\BC^*)-$antilinear involution $K_{G\times \BC^*}(\cN )\to
K_{G\times \BC^*}(\cN ), x\mapsto \bar x$.

Then using usual Kazhdan-Lusztig machinery the basis $\{ C(\lambda )\}$ was
defined. Thus $C(\lambda)$ is a unique selfdual element of the form 
$AJ(\lambda)+\sum_{\mu <\lambda}b_{\mu ,\lambda}AJ(\mu )$ where $b_{\mu ,
\lambda}\in v^{-1}\BZ [v^{-1}]$). The main Conjecture of \cite{O} is that for 
any $\lambda \in X_+$ the support of $C(\lambda)$ is the closure of nilpotent
orbit $\cO_{\lambda}$ and $C(\lambda)|_{\cO_{\lambda}}$ represents up to sign
class of irreducible $G-$equivariant bundle on $\cO_{\lambda}$. In this way
one should recover Lusztig's bijection between dominant weights $X_+$
and pairs consisting of nilpotent orbit and $G-$equivariant irreducible bundle
on it (see \cite{L4}, \cite{B}, \cite{B1}, \cite{V} for various approaches
to Lusztig's bijection). 
All these Conjectures are now known to be true thanks to the work of
R.~Bezrukavnikov.

Now let $e^{\lambda}\in K_G(\cN )$ be the image of $AJ(\lambda )$ under
forgetting map $K_{G\times \BC^*}(\cN )\to K_G(\cN )$. By definition
$e^{\lambda}$ can be constructed as follows: let $\cL (\lambda )$ be the line
bundle on $G/B$ corresponding to the weight $\lambda$ (we choose notations
in such a way that $H^0(G/B, \cL (\lambda ))\ne 0$ for {\em dominant}
$\lambda$); then $e^\lambda =[sp_*\pi^*\cL (\lambda )]$ where $\pi :T^*G/B\to
G/B$ is natural projection and $sp: T^*G/B\to \cN$ is the Springer resolution.
This definition makes sense for any (not necessarily dominant) weight $\lambda$
and one knows that $e^{w\lambda}=e^\lambda$ for any $w\in W, \lambda \in X$.

\subsection{McGovern's formula} \label{Mc}
In this note we are especially interested in weights corresponding to the
trivial bundles on nilpotent orbits under Lusztig's bijection. One of the
Conjectures in \cite{O} states that these weights are exactly $\cL (\cO )$
where $\cL$ is defined in the Introduction. Moreover, $C(\cL (\cO ))$ should
represent a class $j_*(\BC_\cO )$ where $j: \cO \to \cN$ is natural inclusion
and $\BC_\cO$ is trivial $G-$equivariant bundle on $\cO$. Again this is known
to be true thanks to the work of R.~Bezrukavnikov.

Now there is a simple formula for $j_*(\BC_\cO )$ due to W.~McGovern. Namely,
let $h$ be the Dynkin diagram of the orbit $\cO$. Then it defines a grading
of Lie algebra $\fg =\bigoplus_{i\in \BZ}\fg_i$ where $\fg_i=\{ x\in \fg |
[h, x]=ix\}$. Furthermore $\fg_{\ge 0}=\bigoplus_{i\ge 0}\fg_i$ is a parabolic
subalgebra of $\fg$ and $\fg_{\ge 2}=\bigoplus_{i\ge 2}\fg_i$ is a module over
$\fg_{\ge 0}$. Let $G_{\ge 0}$ be a parabolic subgroup with Lie algebra
$\fg_{\ge 0}$ and let $M=G\times_{G_{\ge 0}}\fg_{\ge 2}$ be the homogeneous
bundle on $G/G_{\ge 0}$ corresponding to the $G_0-$module $\fg_{\ge 2}$. There
is a natural map $r: G\times_{G_{\ge 0}}\fg_{\ge 2}\to \fg$. An image of $r$ is
exactly $\bar \cO$ and moreover this map is proper and generically one to one
and therefore is a resolution of singularities of $\cO$, see e.g. \cite{McG}.
In {\em loc. cit.} McGovern proved that $[j_*\BC_{\cO}]=[r_*\BC_M]$.

Now let $R_{+,0}$ (resp. $R_{+,1}$) be the subset of all positive roots such
that $\alpha (h)=0$ (resp. $\alpha (h)=1$). Then the Koszul complex gives the
following formula:
$$
[j_*\BC_{\cO}]=\prod_{\alpha \in R_{+,0}\cup R_{+,1}}(e^0-e^\alpha ) \eqno{(*)}
$$
(see {\em loc. cit.} for details). This implies the following algorithm
for computing $\cL (\cO )$:

1) Multiply the brackets in the right hand side of (*) using usual rule
$e^{\lambda}e^{\mu}=e^{\lambda +\mu}$.

2) In expression from 1) replace each $e^\lambda$ by $e^{w\lambda}$ where
$w\in W$ and $w\lambda$ is dominant. Then make all possible cancelations.
(Warning: Step 1) and Step 2) don't commute!)

3) In expression from 2) find the leading term $\pm e^{\lambda}$ that is such 
term that for any other term $e^\mu$ the inequality $\lambda >\mu$ holds 
(existence of such $\lambda$ is a consequence of the results in \cite{B1}). 
This $\lambda$ is exactly $\cL (\cO )$.

Unfortunately this algorithm is completely impractical for groups of large
rank.

\section{Theorems}
The results of this section are probably well known to experts. Moreover,
Theorem \ref{GL} was stated in \cite{O} without proof.

\subsection{Richardson resolutions} \label{Rich}
It follows from the Theorem of Hinich and
Paniushev \cite{H}, \cite{P} on rationality of singularities of normalizations
of closures of nilpotent orbits that $j_*\BC_{\cO}=r_*\BC_{M}$ where
$r: M\to \bar \cO$ is a resolution of singularities of $\bar \cO$. One obtains
from this McGovern's formula using canonical resolution of a nilpotent orbit.
Now let $P$ be a parabolic subgroup of $G$. The image of the moment
map $m_P: T^*G/P\to \fg^*=\fg$ is the closure of the nilpotent orbit $\cO (P)$
and by the well known Theorem due to R.~W.~Richardson the map $m_P: T^*G/P \to
\overline{\cO (P)}$ is generically finite to one.

Let $G^\vee$ be a Langlands dual group of $G$. By definition there is a
bijection between simple roots for $G$ and $G^\vee$. In particular we can
attach to any parabolic subgroup $P\subset G$ a Levi subgroup $L^\vee_P
\subset G^\vee$. Recall that any coweight for $G^\vee$ is by definition a
weight for $G$.

{\bf Theorem.} {\em Let $P\subset G$ be a parabolic subgroup such that the
map $m_P: T^*G/P \to \overline{\cO (P)}$ is birational. Then $\cL (\cO_P)$
equals to the Dynkin diagram of the principal nilpotent element in $L^\vee_P$.}

{\bf Proof.} By remarks above we know that $j_*\BC_{\cO_P}=m_*\BC_{T^*G/P}$.
Let $R_+(L)\subset R_+$ denote the subset of positive roots of subgroup $L$.
Again the Koszul complex gives a formula
$$
m_*\BC_{T^*G/P}=\prod_{\alpha \in R_+(L)}(1-e^\alpha )
$$
interpreted in the same way as McGovern's formula \ref{Mc}.
Let $\rho_L=\frac{1}{2}\sum_{\alpha \in R_+(L)}\alpha$ and let $W_L\subset W$
denote Weyl group of $L$. The Weyl denominator formula gives
$$
\prod_{\alpha \in R_+(L)}(1-e^\alpha )=\sum_{w\in W_L}det(w)e^{\rho_L-w\rho_L}.
$$
In particular we have a leading term $e^{2\rho_L}$ corresponding to summand
with $w=w_0(L)$ --- the element of largest length in $W_L$. It cannot cancel 
with anything else since the scalar product $\langle 2\rho_L, 2\rho_L \rangle$
is clearly bigger than any scalar product $\langle \rho_L-w\rho_L, \rho_L-
w\rho_L\rangle$ for $w\ne w_0(L)$. By the same reason this term is the unique 
candidate for the leading term since $\lambda >\mu$ for dominant $\lambda$ and
$\mu$ implies $\langle \lambda ,\lambda \rangle >\langle \mu ,\mu \rangle$. 
Since we know that the leading term exists (this is Bezrukavnikov's result as 
we mentionned in \ref{Mc}) it should be equal to $e^{x\rho_L}$ where $x\in W$ 
is such that $x\rho_L$ is dominant.

Finally weight $2\rho_L$ considered as a coweight for $G^\vee$ is clearly the
Dynkin diagram of the regular nilpotent element in $L^\vee$ since the Dynkin
diagram of regular nilpotent element is always the sum of positive coweights,
see \cite{Ko}. \sq

\subsection{} \label{doubly}
We will say that nilpotent orbit is a {\em strongly Richardson orbit} if it
admits a desingularization via momentum map $m_P: T^*G/P\to \bar \cO$. We have
the following

{\bf Theorem.} {\em Let $\cO =\cO_P$ be a strongly Richardson orbit. Suppose
that the $G^\vee -$orbit $\cO'$ of regular nilpotent element of $L^\vee_P$ is 
also strongly Richardson in $G^\vee$. Then $\cL (\cO')=$ Dynkin diagram of 
$\cO$ and $\cL (\cO )=$ Dynkin diagram of $\cO'$.}

{\bf Proof.} The equality $\cL (\cO )=$ Dynkin diagram of $\cO'$ is immediate
from the Theorem \ref{Rich}.

Both nilpotent orbits $\cO$ and $\cO'$ being Richardson orbits
are {\em special}, see \cite{Sp, CM}. N.~Spaltenstein defined an order
reversing involutive bijection $d$ between the sets of special nilpotent
orbits for $G$ and $G^\vee$, see \cite{Sp}. It follows from the definition
\cite{Sp} that $d(\cO')=\cO$. Consequently $d(\cO)=\cO'$ and equality
$\cL (\cO')=$ Dynkin diagram of $\cO$ follows by symmetry. \sq

We don't know the answer to the following

{\bf Question.} Is it true that any Richardson orbit is strongly Richardson?

The nilpotent orbit is called even if its Dynkin diagram is divisible by 2 as
element of coweight lattice, see e.g. \cite{CM}.

{\bf Corollary.} {\em Suppose the orbit $\cO$ is even and that $\cL (\cO )$
is divisible by 2 as element of weight lattice. Then $\cL (\cO)$ is a Dynkin
diagram of some nilpotent orbit $\cO'$ for $G^\vee$ and $\cL (\cO')=$ Dynkin
diagram of $\cO$.}

{\bf Proof.} This is clear since even nilpotent orbits are strongly
Richardson, see e.g. \cite{McG}.\sq

A lot of examples for this Corollary can be found in the Tables in the end of
this note. In fact this is the only general statement we can prove (or even 
just formulate) about these Tables.

\subsection{} It is well known that nilpotent orbits in $GL_n$ are numbered by
partitions of $n$: to the partition $p_1\ge p_2\ge \ldots$ one associates an
orbit of nilpotent matrices with Jordan blocks of size $p_1, p_2, \ldots$.
One finds Dynkin diagram of the orbit $\cO =\cO (p_1, p_2, \ldots )$ as
follows, see e.g. \cite{CM}: order $n$ numbers $p_1-1, p_1-3, \ldots ,
-p_1+1, p_2-1, \ldots ,-p_2+1, \ldots$ in decreasing order, then label of
$i-$th vertice of Dynkin diagram will be difference of $i-$th and $i+1-$th
numbers.

\subsection{} \label{GL1}
The following Lemma is well known.

{\bf Lemma.} {\em For $G=GL_n$ and a parabolic subgroup $P\subset G$ the map
$m_P: T^*G/P \to \overline{\cO (P)}$ is birational.}

{\bf Proof.} The isotropy group of any nilpotent element in $GL_n$ is
connected. \sq

\subsection{} \label{GL}
Let $\cO =\cO (p_1, p_2, \ldots )$ be a nilpotent orbit. Let
$p_1', p_2', \ldots$ be a partition dual to $p_1, p_2, \ldots$ and let
$\cO'=\cO (p_1', p_2', \ldots )$ be the orbit dual to $\cO$.

{\bf Theorem.} {\em The weight $\cL (\cO )$ equals to the Dynkin diagram
of the dual orbit $\cO'$ considered as a weight for $GL_n$.}

{\bf Proof.} It is well known that any nilpotent orbit $\cO$ in $G=GL_n$ is a
Richardson orbit and it follows from the Lemma \ref{GL1} that it is strongly
Richardson. It follows from the proof of the Theorem \ref{doubly} that
$\cL (\cO)$ is Dynkin diagram of $d(\cO)$ (here $d$ is Spaltenstein duality).
Now the result follows from description of $d$ for the group $GL_n$ in 
\cite{Sp}. \sq

\subsection{Remarks} (i) As we mentioned in the Introduction the calculation of
the map $\cL$ is equivalent to the calculation of canonical distinguished
involutions in affine Weyl group. For type $A_n$ canonical distinguished
involutions were recently calculated by Nanhua Xi \cite{Xi} by a completely
different method. Surely his answer coincides with ours, but he did not mention
a relation with Dynkin diagrams of nilpotent orbits.

(ii) For group $G$ of type different from $A_n$ it {\em is not} true in general
that weight attached to a nilpotent orbit and considered as coweight for
Langlands dual group $G^\vee$ is a Dynkin diagram for some nilpotent orbit of
$G^\vee$, see tables in the end of this note. But calculations made by
Pramod Achar (private communinication) suggest some evidence for the positive
answer to the following

{\bf Question.} Is it true that any Dynkin diagram considered as a weight for
the dual group corresponds under Lusztig's bijection to the local system on 
some nilpotent orbit? (Here local system is a coherent equivariant sheaf on 
nilpotent orbit such that the corresponding representation of isotropy group 
factors through finite quotient.)

\section{Tables}
The Tables below contain results of our calculations of map $\cL$ for groups
of small rank. We have almost complete results for groups of rank $\le 7$
(there are some gaps for groups of types $B_6$, $C_6$, $B_7$, $C_7$ and $E_7$)
and partial results in rank 8 (for groups $D_8$ and $E_8$). The Tables are
organized as follows: the tables for classical groups consist of 4 columns,
in the first column we give the partition identifying the nilpotent orbit 
$\cO$ (see e.g. \cite{CM}), in the second column we give the Dynkin diagram of
$\cO$, in the third column the weight $\cL(\cO)$ is contained and last column 
contains the square of the length of $\cL(\cO)$ (we normalize scalar product
by $(\alpha,\alpha)=2$ for a short root $\alpha$); for exceptional groups the 
tables consist of 3 columns, the first column contains Dynkin diagram of the 
nilpotent orbit $\cO$, the second column contains $\cL(\cO)$ and the third 
column contains the square of the length of $\cL(\cO)$. In a few cases we were
unable to calculate $\cL(\cO)$ but we could predict the value $\cL(\cO)$ using
duality reasoning, such cases are marked by a question mark. 

\newpage

\hspace{50pt}
\begin{minipage}[t]{2.4in}
\centerline{\large\bf so$_5$}
\vspace{7pt}                              
$
\begin{array}{c|c|c|c}
\mbox{Partition} & \mbox{Diagram} & \mbox{Weight} & \\ \hline
(5)     & \bb22 & \bb00 & 0 \\
(3,1^2) & \bb20 & \bb10 & 2 \\
(2^2,1) & \bb01 & \bb12 & 10 \\
(1^5)   & \bb00 & \bb22 & 20
\end{array}
$

\vspace{0.5cm}
\centerline{\large\bf so$_7$}
\vspace{7pt}
$
\begin{array}{c|c|c|c}
\mbox{Partition} & \mbox{Diagram} & \mbox{Weight} & \\ \hline
(7)       & \bbb222 & \bbb000 & 0  \\
(5,1^2)   & \bbb220 & \bbb100 & 2  \\
(3^2,1)   & \bbb020 & \bbb002 & 6  \\
(3,2^2)   & \bbb101 & \bbb020 & 16  \\
(3,1^4)   & \bbb200 & \bbb210 & 20 \\
(2^2,1^3) & \bbb010 & \bbb112 & 28  \\
(1^7)     & \bbb000 & \bbb222 & 70
\end{array}
$

\vspace{0.5cm}
\centerline{\large\bf so$_9$}
\vspace{7pt}
$
\begin{array}{c|c|c|c}
\mbox{Partition} & \mbox{Diagram} & \mbox{Weight} & \\ \hline
(9)         & \bbbb2222 & \bbbb0000 & 0  \\
(7,1^2)     & \bbbb2220 & \bbbb1000 & 2  \\
(5,3,1)     & \bbbb2020 & \bbbb0010 & 6  \\
(4^2,1)     & \bbbb0201 & \bbbb1002 & 14 \\
(5,2^2)     & \bbbb2101 & \bbbb0200 & 16 \\
(3^3)       & \bbbb0020 & \bbbb0110 & 18 \\
(5,1^4)     & \bbbb2200 & \bbbb2100 & 20 \\
(3^2,1^3)   & \bbbb0200 & \bbbb2002 & 24 \\
(3,2^2,1^2) & \bbbb0101 & \bbbb0202 & 40 \\
(2^4,1)     & \bbbb0001 & \bbbb1112 & 60 \\
(3,1^6)     & \bbbb2000 & \bbbb2210 & 70 \\
(2^2,1^5)   & \bbbb0100 & \bbbb2112 & 78 \\
(1^9)       & \bbbb0000 & \bbbb2222 & 168 
\end{array}
$
\end{minipage}
\hspace{40pt}
\begin{minipage}[t]{2.3in}
\centerline{\large\bf sp$_4$}
\vspace{7pt}
$
\begin{array}{c|c|c|c}
\mbox{Partition} & \mbox{Diagram} & \mbox{Weight} & \\ \hline
(4)     & \s22 & \s00 & 0  \\
(2^2)   & \s02 & \s01 & 2  \\
(2,1^2) & \s10 & \s21 & 10 \\
(1^4)   & \s00 & \s22 & 20
\end{array}
$

\vspace{0.5cm}
\centerline{\large\bf sp$_6$}
\vspace{7pt}
$
\begin{array}{c|c|c|c}
\mbox{Partition} & \mbox{Diagram} & \mbox{Weight} & \\ \hline
(6)       & \sp222 & \sp000 & 0  \\
(4,2)     & \sp202 & \sp010 & 2  \\
(3^2)     & \sp020 & \sp101 & 6  \\
(2^3)     & \sp002 & \sp020 & 8  \\
(4,1^2)   & \sp210 & \sp210 & 10 \\
(2^2,1^2) & \sp010 & \sp220 & 20 \\
(2,1^4)   & \sp100 & \sp230 & 34 \\
(1^6)     & \sp000 & \sp222 & 56
\end{array}
$

\vspace{0.5cm}
\centerline{\large\bf sp$_8$}
\vspace{7pt}
$
\begin{array}{c|c|c|c}
\mbox{Partition} & \mbox{Diagram} & \mbox{Weight} & \\ \hline
(8)       & \spp2222 & \spp0000 & 0  \\
(6,2)     & \spp2202 & \spp0100 & 2  \\
(4^2)     & \spp0202 & \spp0001 & 4  \\
(4,2^2)   & \spp2002 & \spp0200 & 8  \\
(6,1^2)   & \spp2210 & \spp2100 & 10 \\
(3^2,2)   & \spp0110 & \spp0020 & 12 \\
(4,2,1^2) & \spp2010 & \spp2200 & 20 \\
(2^4)     & \spp0002 & \spp0201 & 20 \\
(3^2,1^2) & \spp0200 & \spp2101 & 22 \\
(4,1^4)   & \spp2100 & \spp2300 & 34 \\
(2^3,1^2) & \spp0010 & \spp2201 & 36 \\
(2^2,1^4) & \spp0100 & \spp2301 & 54 
\end{array}
$
\end{minipage}

\newpage

\hspace{30pt}
\begin{minipage}[t]{2.8in}
\centerline{\large\bf{so$_{11}$}} 
\vspace{7pt}
$
\begin{array}{c|c|c|c}
\mbox{Partition} & \mbox{Diagram} & \mbox{Weight} & \\ \hline
(11)        & \bbbbb22222 & \bbbbb00000 & 0   \\
(9,1^2)     & \bbbbb22220 & \bbbbb10000 & 2   \\
(7,3,1)     & \bbbbb22020 & \bbbbb00100 & 6   \\
(5^2,1)     & \bbbbb02020 & \bbbbb00002 & 10   \\
(7,2^2)     & \bbbbb22101 & \bbbbb02000 & 16   \\
(5,3^2)     & \bbbbb20020 & \bbbbb01100 & 18   \\
(7,1^4)     & \bbbbb22200 & \bbbbb21000 & 20  \\
(5,3,1^3)   & \bbbbb20200 & \bbbbb20010 & 24  \\
(4^2,3)     & \bbbbb01101 & \bbbbb00110 & 26  \\
(4^2,1^3)   & \bbbbb02010 & \bbbbb11002 & 32  \\
(3^3,1^2)   & \bbbbb00200 & \bbbbb10110 & 36  \\
(5,2^2,1^2) & \bbbbb21010 & \bbbbb02010 & 40  \\
(3^2,2^2,1) & \bbbbb01010 & \bbbbb11110 & 60  \\
(5,1^6)     & \bbbbb22000 & \bbbbb22100 & 70  \\
(3^2,1^5)   & \bbbbb02000 & \bbbbb22002 & 74  \\
(3,2^4)     & \bbbbb10001 & \bbbbb02020 & 80  \\
(3,2^2,1^4) & \bbbbb10100 & \bbbbb20202 & 90  \\
(2^4,1^3)   & \bbbbb00010 & \bbbbb11112 & 110  \\
(3,1^8)     & \bbbbb20000 & \bbbbb22210 & 168  \\
(2^2,1^7)   & \bbbbb01000 & \bbbbb22112 & 176  \\
(1^{11})    & \bbbbb00000 & \bbbbb22222 & 330
\end{array}
$

\vspace{0.5cm}
\centerline{\large\bf{so$_{13}$}} 
\vspace{7pt}
$
\begin{array}{c|c|c|c}
\mbox{Partition} & \mbox{Diagram} & \mbox{Weight} & \\ \hline
(13)          & \bbbbbb222222 & \bbbbbb000000 & 0  \\
(11,1^2)      & \bbbbbb222220 & \bbbbbb100000 & 2  \\
(9,3,1)       & \bbbbbb222020 & \bbbbbb001000 & 6  \\
(7,5,1)       & \bbbbbb202020 & \bbbbbb000010 & 10 \\ 
(5,7,1)       & \bbbbbb202020 & \bbbbbb000010 & 10 \\
(9,2^2)       & \bbbbbb222101 & \bbbbbb020000 & 16  
\end{array}
$
\end{minipage}
\hspace{20pt}
\begin{minipage}[t]{2.8in}
\centerline{\large\bf{sp$_{8}$}} 
\vspace{7pt}
$
\begin{array}{c|c|c|c}
\mbox{Partition} & \mbox{Diagram} & \mbox{Weight} & \\ \hline
(2,1^6)   & \spp1000 & \spp2221 & 84 \\
(1^8)     & \spp0000 & \spp2222 & 120
\end{array}
$

\vspace{0.8cm}
\centerline{\large\bf{sp$_{10}$}} 
\vspace{7pt}
$
\begin{array}{c|c|c|c}
\mbox{Partition} & \mbox{Diagram} & \mbox{Weight} & \\ \hline
(10)        & \sppp22222 & \sppp00000 & 0   \\
(8,2)       & \sppp22202 & \sppp01000 & 2   \\
(6,4)       & \sppp20202 & \sppp00010 & 4   \\
(6,2^2)     & \sppp22002 & \sppp02000 & 8   \\
(5^2)       & \sppp02020 & \sppp10001 & 8   \\
(4^2,2)     & \sppp02002 & \sppp01010 & 10  \\
(8,1^2)     & \sppp22210 & \sppp21000 & 10  \\
(4,3^2)     & \sppp10110 & \sppp10110 & 18  \\
(6,2,1^2)   & \sppp22010 & \sppp22000 & 20  \\
(4,2^3)     & \sppp20002 & \sppp02010 & 20  \\
(4^2,1^2)   & \sppp02010 & \sppp21010 & 22  \\
(3^2,2^2)   & \sppp01010 & \sppp01101 & 24  \\
(6,1^4)     & \sppp22100 & \sppp23000 & 34  \\
(4,2^2,1^2) & \sppp20010 & \sppp22010 & 36  \\
(3^2,2,1^2) & \sppp01100 & \sppp21101 & 40  \\
(2^5)       & \sppp00002 & \sppp02020 & 40  \\
(4,2,1^4)   & \sppp20100 & \sppp23010 & 54  \\
(3^2,1^4)   & \sppp02000 & \sppp22101 & 58  \\
(2^4,1^2)   & \sppp00010 & \sppp22020 & 60  \\ 
(2^3,1^4)   & \sppp00100 & \sppp23020 & 82  \\
(4,1^6)     & \sppp21000 & \sppp22210 & 84  \\
(2^2,1^6)   & \sppp01000 & \sppp22220 & 120 \\
(2,1^8)     & \sppp10000 & \sppp22230 & 164 \\
(1^{10})    & \sppp00000 & \sppp22222 & 220
\end{array}
$
\end{minipage}

\newpage

\hspace{10pt}
\begin{minipage}[t]{8cm}
\centerline{\large\bf{so$_{13}$}}
\vspace{7pt}
$
\begin{array}{c|c|c|c}
\mbox{Partition} & \mbox{Diagram} & \mbox{Weight} & \\ \hline
(7,3^2)       & \bbbbbb220020 & \bbbbbb011000 & 18 \\
(6^2,1)       & \bbbbbb020201 & \bbbbbb100002 & 18   \\
(9,1^4)       & \bbbbbb222200 & \bbbbbb210000 & 20 \\
(5^2,3)       & \bbbbbb020020 & \bbbbbb010010 & 22  \\
(7,3,1^3)     & \bbbbbb220200 & \bbbbbb200100 & 24 \\
(5^2,1^3)     & \bbbbbb020200 & \bbbbbb200002 & 24  \\
(5,4^2)       & \bbbbbb101101 & \bbbbbb000200 & 32  \\
(5,3^2,1^2)   & \bbbbbb200200 & \bbbbbb101100 & 36  \\
(7,2^2,1)     & \bbbbbb221010 & \bbbbbb020100 & 40 \\
(4^2,3,1^2)   & \bbbbbb011010 & \bbbbbb100110 & 44  \\
(5,3,2^2,1)   & \bbbbbb201010 & \bbbbbb111100 & 60  \\
(3^4,1)       & \bbbbbb000200 & \bbbbbb002002 & 60  \\
(4,2^2,1)     & \bbbbbb020001 & \bbbbbb111002 & 64  \\
(7,1^6)       & \bbbbbb222000 & \bbbbbb221000 & 70 \\
(5,3,1^5)     & \bbbbbb202000 & \bbbbbb220010 & 74  \\
(5,2^4)       & \bbbbbb210001 & \bbbbbb020200 & 80  \\
(4^2,1^5)     & \bbbbbb020100 & \bbbbbb211002 & 82  \\
(3^3,2^2)     & \bbbbbb001010 & \bbbbbb020110 & 82  \\
(3^3,1^4)     & \bbbbbb002000 & \bbbbbb210110 & 86  \\ 
(5,2^2,1^4)   & \bbbbbb210100 & \bbbbbb202010 & 90  \\
(3^2,2^2,1^3) & \bbbbbb010100 & \bbbbbb111110 & 110  \\
(3,2^4,1^2)   & \bbbbbb100010 & \bbbbbb020210 & 138  \\
(5,1^8)       & \bbbbbb220000 & \bbbbbb222100 & 168  \\
(3^2,1^7)     & \bbbbbb020000 & \bbbbbb222002 & 172  \\
(2^6,1)       & \bbbbbb000001 & \bbbbbb111112 & 182  \\
(3,1^{10})    & \bbbbbb200000 & \bbbbbb222210 & 330 \\
(3,2^2,1^6)   & \bbbbbb101000 & &           \\
(2^4,1^5)     & \bbbbbb000100 & & 
\end{array}
$ 
\end{minipage}
\hspace{20pt}
\begin{minipage}[t]{8cm}
\centerline{\large\bf{sp$_{12}$}} 
\vspace{7pt}
$
\begin{array}{c|c|c|c}
\mbox{Partition} & \mbox{Diagram} & \mbox{Weight} &  \\ \hline
(12)          & \ssp222222 & \ssp000000 & 0  \\
(10,2)        & \ssp222202 & \ssp010000 & 2  \\ 
(8,4)         & \ssp220202 & \ssp000100 & 4   \\
(6^2)         & \ssp020202 & \ssp000001 & 6   \\
(8,2^2)       & \ssp222002 & \ssp020000 & 8   \\
(10,1^2)      & \ssp222210 & \ssp210000 & 10  \\
(6,4,2)       & \ssp202002 & \ssp010100 & 10  \\
(5^2,2)       & \ssp020110 & \ssp001010 & 14  \\
(4^3)         & \ssp002002 & \ssp000200 & 16  \\
(6,3^2)       & \ssp210110 & \ssp101100 & 18  \\
(6,2^3)       & \ssp220002 & \ssp020100 & 20  \\
(8,2,1^2)     & \ssp222101 & \ssp220000 & 20  \\
(6,4,1^2)     & \ssp202010 & \ssp210100 & 22  \\
(4^2,2^2)     & \ssp020002 & \ssp020001 & 22  \\
(5^2,1^2)     & \ssp020200 & \ssp210001 & 24  \\
(8,1^4)       & \ssp222100 & \ssp230000 & 34  \\
(6,2^2,1^2)   & \ssp220010 & \ssp220100 & 36  \\
(4,2^2)       & \ssp200002 & \ssp020200 & 40  \\
(3^4)         & \ssp000200 & \ssp101101 & 40  \\
(4,3^2,1^2)   & \ssp101100 & \ssp210200 & 42  \\
(3^2,2^3)     & \ssp010010 & \ssp020020 & 44  \\
(6,2,1^4)     & \ssp220100 & \ssp230100 & 54  \\
(4^2,1^4)     & \ssp020100 & \ssp230001 & 56  \\
(4,3^2,2)     & \ssp101010 & \ssp002002 & 60  \\
(4,2^3,1^2)   & \ssp200010 & \ssp220200 & 60  \\
(3^2,2^2,1^2) & \ssp010100 & \ssp220020 & 64  \\
(2^6)         & \ssp000002 & \ssp020201 & 70  \\
(4,2^2,1^4)   & \ssp200100 & \ssp230200 & 82  
\end{array}
$
\end{minipage}

\newpage

\begin{minipage}[t]{8.5cm}
\centerline{\large\bf{so$_{13}$}}
\vspace{7pt}
$
\begin{array}{c|c|c|c}
\mbox{Partition} & \mbox{Diagram} & \mbox{Weight} & \\ \hline
(2^2,1^9)     & \bbbbbb010000 & &  \\
(1^{13})      & \bbbbbb000000 & \bbbbbb222222 & 572
\end{array}
$

\vspace{0.5cm}
\centerline{\large\bf so$_{15}$}
\vspace{7pt}
$
\begin{array}{c|c|c|c}
\mbox{Partition} & \mbox{Diagram} & \mbox{Weight} & \\ \hline
(15)          & \bbbbbbb2222222 & \bbbbbbb0000000 & 0  \\
(13,1^2)      & \bbbbbbb2222220 & \bbbbbbb1000000 & 2  \\
(11,3,1)      & \bbbbbbb2222020 & \bbbbbbb0010000 & 6  \\
(9,5,1)       & \bbbbbbb2202020 & \bbbbbbb0000100 & 10  \\
(7^2,1)       & \bbbbbbb0202020 & \bbbbbbb0000002 & 14  \\
(11,2^2)      & \bbbbbbb2222101 & \bbbbbbb0200000 & 16  \\
(9,3^2)       & \bbbbbbb2220020 & \bbbbbbb0110000 & 18  \\
(11,1^4)      & \bbbbbbb2222200 & \bbbbbbb2100000 & 20 \\
(7,5,3)       & \bbbbbbb2020020 & \bbbbbbb0100100 & 22 \\
(7,5,1^3)     & \bbbbbbb2020200 & \bbbbbbb2000010 & 28 \\
(6^2,3)       & \bbbbbbb0201101 & \bbbbbbb0010010 & 30 \\
(7,4^2)       & \bbbbbbb2101101 & \bbbbbbb0002000 & 32 \\
(5^3)         & \bbbbbbb0020020 & \bbbbbbb0001100 & 34 \\
(7,3^2,1^2)   & \bbbbbbb2200200 & \bbbbbbb1011000 & 36 \\
(6^2,1^3)     & \bbbbbbb0202010 & \bbbbbbb1100002 & 36 \\
(9,2^2,1^2)   & \bbbbbbb2221010 & \bbbbbbb0201000 & 40 \\
(5^2,3,1^2)   & \bbbbbbb0200200 & \bbbbbbb1010010 & 40 \\
(5,4^2,1^2)   & \bbbbbbb1011010 & \bbbbbbb0101010 & 56 \\
(5,3^3,1)     & \bbbbbbb2000200 & \bbbbbbb0020010 & 60 \\
(7,3,2^2,1)   & \bbbbbbb2201010 & \bbbbbbb1111000 & 60 \\
(5^2,2^2,1)   & \bbbbbbb0201010 & \bbbbbbb1110010 & 64 \\
(4^2,3^2,1)   & \bbbbbbb0101010 & \bbbbbbb0011002 & 68  \\
(9,1^6)       & \bbbbbbb2222000 & \bbbbbbb2210000 & 70 \\
(7,3,1^5)     & \bbbbbbb2202000 & \bbbbbbb2200100 & 74 \\
(5^2,1^5)     & \bbbbbbb0202000 & \bbbbbbb2200002 & 78 
\end{array}
$
\end{minipage}
\hspace{10pt}
\begin{minipage}[t]{8.5cm}
\centerline{\large\bf{sp$_{12}$}} 
\vspace{7pt}
$
\begin{array}{c|c|c|c}
\mbox{Partition} & \mbox{Diagram} & \mbox{Weight} &  \\ \hline
(6,1^6)       & \ssp221000 & \ssp222100 & 84  \\
(3^2,2,1^4)   & \ssp011000 & \ssp230020 & 86  \\
(2^5,1^2)     & \ssp000010 & \ssp220201 & 94  \\
(4,2,1^6)     & \ssp201000 & \ssp222200 & 120 \\
(3^2,1^6)     & \ssp020000 & \ssp222101 & 122 \\
(4,1^8)       & \ssp210000 & \ssp222300 & 164 \\
(2^4,1^4)     & \ssp000100 & &           \\
(2^3,1^6)     & \ssp001000 & &           \\
(2^2,1^8)     & \ssp010000 & &           \\
(2,1^{10})    & \ssp100000 & &           \\
(1^{12})      & \ssp000000 & \ssp222222 & 364          
\end{array}
$

\vspace{0.5cm}
\centerline{\large\bf{sp$_{14}$}} 
\vspace{7pt}
$
\begin{array}{c|c|c|c}
\mbox{Partition} & \mbox{Diagram} & \mbox{Weight} & \\ \hline
(14)          & \sssp2222222 & \sssp0000000 & 0  \\
(12,2)        & \sssp2222202 & \sssp0100000 & 2  \\
(10,4)        & \sssp2220202 & \sssp0001000 & 4  \\
(8,6)         & \sssp2020202 & \sssp0000010 & 6  \\
(10,2^2)      & \sssp2222002 & \sssp0200000 & 8  \\
(12,1^2)      & \sssp2222210 & \sssp2100000 & 10 \\
(8,4,2)       & \sssp2202002 & \sssp0101000 & 10 \\
(7^2)         & \sssp0202020 & \sssp1000001 & 10 \\
(6^2,2)       & \sssp0202002 & \sssp0100010 & 12 \\
(6,4^2)       & \sssp2002002 & \sssp0002000 & 16 \\
(8,3^2)       & \sssp2210110 & \sssp1011000 & 18 \\
(10,2,1^2)    & \sssp2222010 & \sssp2200000 & 20 \\
(8,2^3)       & \sssp2220002 & \sssp0201000 & 20 \\
(5^2,4)       & \sssp0110110 & \sssp0000200 & 20  \\
(8,4,1^2)     & \sssp2202010 & \sssp2101000 & 22 \\
(6,4,2^2)     & \sssp2020002 & \sssp0200010 & 22 
\end{array}
$
\end{minipage}

\newpage

\begin{minipage}[t]{8.5cm}
\centerline{\large\bf so$_{15}$}
\vspace{7pt}
$
\begin{array}{c|c|c|c}
\mbox{Partition} & \mbox{Diagram} & \mbox{Weight} & \\ \hline
(7,2^4)       & \bbbbbbb2210001 & \bbbbbbb0202000 & 80 \\
(5,3^2,2^2)   & \bbbbbbb2001010 & \bbbbbbb0201100 & 82 \\
(5,3^2,1^4)   & \bbbbbbb2002000 & \bbbbbbb2101100 & 86 \\
(4^2,3,2^2)   & \bbbbbbb0110001 & \bbbbbbb0200110 & 90  \\
(7,2^2,1^4)   & \bbbbbbb2210100 & \bbbbbbb2020100 & 90 \\
(4^2,3,1^4)   & \bbbbbbb0110100 & \bbbbbbb2100110 & 94  \\
(3^5)         & \bbbbbbb0000200 & \bbbbbbb0110110 & 100  \\
(3^4,1^3)     & \bbbbbbb0002000 & \bbbbbbb2002002 & 110  \\
(5,3,2^2,1^3) & \bbbbbbb2010100 & \bbbbbbb1111100 & 110 \\
(4^2,2^2,1^3) & \bbbbbbb0200010 & \bbbbbbb1111002 & 114  \\
(5,2^4,1^2)   & \bbbbbbb2100010 & \bbbbbbb0202100 & 138  \\
(7,1^8)       & \bbbbbbb2220000 & \bbbbbbb2221000 & 168 \\
(5,3,1^7)     & \bbbbbbb2020000 & \bbbbbbb2220010 & 172 \\
(4^2,1^7)     & \bbbbbbb0201000 & \bbbbbbb2211002 & 180  \\
(3^3,1^6)     & \bbbbbbb0020000 & \bbbbbbb2210110 & 184  \\
(5,1^{10})    & \bbbbbbb2200000 & \bbbbbbb2222100 & 330 \\
(3^2,1^9)     & \bbbbbbb0200000 & \bbbbbbb2222002 & 334 \\
(3,1^{12})    & \bbbbbbb2000000 & \bbbbbbb2222210 & 572 \\
(5,2^2,1^6)   & \bbbbbbb2101000 &  & \\
(3^3,2^2,1^2) & \bbbbbbb0010100 & & \\
(3^2,2^4,1)   & \bbbbbbb0100010 & & \\
(3^2,2^2,1^5) & \bbbbbbb0101000 & & \\
(3,2^6)       & \bbbbbbb1000001 & & \\
(3,2^4,1^4)   & \bbbbbbb1000100 & & \\
(3,2^2,1^8)   & \bbbbbbb1010000 & & \\
(2^6,1^3)     & \bbbbbbb0000010 & & \\
(2^4,1^7)     & \bbbbbbb0001000 & & \\
(2^2,1^{11})  & \bbbbbbb0100000 & & \\
(1^{15})      & \bbbbbbb0000000 & \bbbbbbb2222222 & 910 
\end{array}
$
\end{minipage}
\hspace{10pt}
\begin{minipage}[t]{8.5cm}
\centerline{\large\bf sp$_{14}$}
\vspace{7pt}
$
\begin{array}{c|c|c|c}
\mbox{Partition} & \mbox{Diagram} & \mbox{Weight} & \\ \hline
(6^2,1^2)     & \sssp0202010 & \sssp2100010 & 24 \\
(5^2,2^2)     & \sssp0201010 & \sssp0110001 & 26  \\
(4^3,2)       & \sssp0020002 & \sssp0101010 & 28  \\  
(6,3^2,2)     & \sssp2101010 & \sssp0020010 & 30 \\
(10,1^4)      & \sssp2222100 & \sssp2300000 & 34 \\
(8,2^2,1^2)   & \sssp2220010 & \sssp2201000 & 36 \\
(6,4,2,1^2)   & \sssp2020010 & \sssp2200010 & 38 \\
(6,2^4)       & \sssp2200002 & \sssp0202000 & 40 \\
(4^2,3^2)     & \sssp0101010 & \sssp1011010 & 40  \\
(4^2,2^3)     & \sssp0200002 & \sssp0201010 & 42  \\
(6,3^2,1^2)   & \sssp2101100 & \sssp2102000 & 42 \\
(5^2,2,1^2)   & \sssp0201100 & \sssp2110001 & 42  \\
(4^3,1^2)     & \sssp0020010 & \sssp2101010 & 44  \\
(4,3^2,2^2)   & \sssp1010010 & \sssp0110110 & 50  \\
(8,2,1^4)     & \sssp2220100 & \sssp2301000 & 54 \\
(6,4,1^4)     & \sssp2020100 & \sssp2300010 & 56 \\
(6,2^3,1^2)   & \sssp2200010 & \sssp2202000 & 60 \\
(5^2,1^4)     & \sssp0202000 & \sssp2210001 & 60  \\ 
(3^4,2)       & \sssp0001100 & \sssp0020020 & 60  \\
(4^2,2^2,1^2) & \sssp0200010 & \sssp2201010 & 62  \\
(4,3^2,2,1^2) & \sssp1010100 & \sssp2110110 & 70  \\
(4,2^5)       & \sssp2000002 & \sssp0202010 & 70  \\ 
(3^4,1^2)     & \sssp0002000 & \sssp2101101 & 76  \\
(6,2^2,1^4)   & \sssp2200100 & \sssp2302000 & 82 \\
(4^2,2,1^4)   & \sssp0200100 & \sssp2301010 & 84  \\
(8,1^6)       & \sssp2221000 & \sssp2221000 & 84 \\
(4,3^2,1^4)   & \sssp1011000 & \sssp2210110 & 92  \\ 
(4,2^4,1^2)   & \sssp2000010 & \sssp2202010 & 94 \\  
(6,2,1^6)     & \sssp2201000 & \sssp2222000 & 120 
\end{array}
$
\end{minipage}

\newpage

\hspace{8.5cm}
\begin{minipage}[t]{8.5cm}
\centerline{\large\bf sp$_{14}$}
\vspace{7pt}
$
\begin{array}{c|c|c|c}
\mbox{Partition} & \mbox{Diagram} & \mbox{Weight} & \\ \hline
(2^7)         & \sssp0000002 & \sssp0202020 & 112 \\
(4^2,1^6)     & \sssp0201000 & \sssp2221010 & 122 \\
(6,1^8)       & \sssp2210000 & \sssp2223000 & 164 \\
(3^2,1^8)     & \sssp0200000 & \sssp2222101 & 222 \\
(4,2^3,1^4)   & \sssp2000100 & & \\
(4,2^2,1^6)   & \sssp2001000 & & \\
(4,2,1^8)     & \sssp2010000 & & \\
(4,1^10)      & \sssp2100000 & & \\
(3^2,2^4)     & \sssp0100010 & & \\
(3^2,2^3,1^2) & \sssp0100100 & & \\
(3^2,2^2,1^4) & \sssp0101000 & & \\
(3^2,2,1^6)   & \sssp0110000 & & \\
(2^6,1^2)     & \sssp0000010 & & \\
(2^5,1^4)     & \sssp0000100 & & \\
(2^4,1^6)     & \sssp0001000 & & \\
(2^3,1^8)     & \sssp0010000 & & \\
(2^2,1^{10})  & \sssp0100000 & & \\
(2,1^{12})    & \sssp1000000 & & \\
(1^{14})      & \sssp0000000 & \sssp2222222 & 560
\end{array}
$
\end{minipage}

\newpage

\hspace{50pt}
\begin{minipage}[t]{6cm}
\centerline{\large\bf so$_6$} 
\vspace{7pt}
$
\begin{array}{c|c|c|c}
\mbox{Partition} & \mbox{Diagram} & \mbox{Weight} &  \\ \hline
(5,1)     & \so222 & \so000 & 0 \\
(3^2)     & \so022 & \so011 & 2 \\
(3,1^3)   & \so200 & \so200 & 4 \\
(2^2,1^2) & \so011 & \so022 & 8 \\
(1^6)     & \so000 & \so222 & 20
\end{array}
$

\vspace{0.5cm}
\centerline{\large\bf so$_8$}
\vspace{7pt}
$
\begin{array}{c|c|c|c}
\mbox{Partition} & \mbox{Diagram} & \mbox{Weight} &  \\ \hline
(7,1)     & \soo2222 & \soo0000 & 0  \\
(5,3)     & \soo2022 & \soo0100 & 2  \\
(5,1^3)   & \soo2200 & \soo2000 & 4  \\
(4^2)_1   & \soo0202 & \soo0002 & 4  \\
(4^2)_2   & \soo0220 & \soo0020 & 4  \\
(3^2,1^2) & \soo0200 & \soo1011 & 6  \\
(3,2^2,1) & \soo1011 & \soo1111 & 14 \\
(3,1^5)   & \soo2000 & \soo2200 & 20 \\
(2^4)_1   & \soo0002 & \soo0202 & 20 \\
(2^4)_2   & \soo0020 & \soo0220 & 20 \\
(2^2,1^4) & \soo0100 & \soo2022 & 24 \\
(1^8)     & \soo0000 & \soo2222 & 56
\end{array}
$

\vspace{0.5cm}
\centerline{\large\bf so$_{10}$}
\vspace{7pt}  
$
\begin{array}{c|c|c|c}
\mbox{Partition} & \mbox{Diagram} & \mbox{Weight} &  \\ \hline
(9,1)       & \sooo22222 & \sooo00000 & 0  \\
(7,3)       & \sooo22022 & \sooo01000 & 2  \\
(7,1^3)     & \sooo22200 & \sooo20000 & 4  
\end{array}
$
\end{minipage}
\hspace{50pt}
\begin{minipage}[t]{7cm}
\centerline{\large\bf so$_{10}$}
\vspace{7pt}
$
\begin{array}{c|c|c|c}
\mbox{Partition} & \mbox{Diagram} & \mbox{Weight} &  \\ \hline
(5^2)       & \sooo02022 & \sooo00011 & 4  \\
(5,3,1^2)   & \sooo20200 & \sooo10100 & 6  \\
(4^2,1^2)   & \sooo02011 & \sooo01011 & 10 \\
(3^3,1)     & \sooo00200 & \sooo00200 & 12 \\
(5,2^2,1)   & \sooo21011 & \sooo11100 & 14 \\
(5,1^5)     & \sooo22000 & \sooo22000 & 20 \\
(3^2,2^2)   & \sooo01011 & \sooo02011 & 20 \\
(3^2,1^4)   & \sooo02000 & \sooo21011 & 22 \\
(3,2^2,1^3) & \sooo10100 & \sooo11111 & 30 \\
(2^4,1^2)   & \sooo00011 & \sooo02022 & 40 \\
(3,1^7)     & \sooo20000 & \sooo22200 & 56 \\
(2^2,1^6)   & \sooo01000 & \sooo22022 & 60 \\
(1^{10})    & \sooo00000 & \sooo22222 & 120
\end{array}
$

\vspace{0.45cm}
\centerline{\large\bf so$_{12}$}
\vspace{7pt}
$
\begin{array}{c|c|c|c}
\mbox{Partition} & \mbox{Diagram} & \mbox{Weight} &  \\ \hline
(11,1)        & \sso222222 & \sso000000 & 0  \\
(9,3)         & \sso222022 & \sso010000 & 2  \\
(9,1^3)       & \sso222200 & \sso200000 & 4  \\
(7,5)         & \sso202022 & \sso000100 & 4  \\
(7,3,1^2)     & \sso220200 & \sso101000 & 6  \\
(6^2)_1       & \sso020220 & \sso000002 & 6  \\ 
(6^2)_2       & \sso020202 & \sso000020 & 6  \\
(5^2,1^2)     & \sso020200 & \sso100011 & 8  \\
(5,3^2,1)     & \sso200200 & \sso002000 & 12 
\end{array}
$
\end{minipage}
        
\newpage
\hspace{7pt}
\vspace{10pt}
\begin{minipage}[t]{8cm}
\centerline {\large\bf so$_{12}$}
\vspace{10pt}
$
\begin{array}{c|c|c|c}
\mbox{Partition} & \mbox{Diagram} & \mbox{Weight} &  \\ \hline
(7,2^2,1)     & \sso221011 & \sso111000 & 14 \\
(4^2,3,1)     & \sso011011 & \sso000200 & 16  \\
(7,1^5)       & \sso222000 & \sso220000 & 20 \\
(5,3,2^2)     & \sso201011 & \sso020100 & 20 \\
(4^2,2^2)_1   & \sso020020 & \sso020002 & 22  \\
(4^2,2^2)_2   & \sso020002 & \sso020020 & 22  \\
(5,3,1^4)     & \sso202000 & \sso210100 & 22 \\
(3^4)         & \sso000200 & \sso011011 & 24  \\
(4^2,1^4)     & \sso020100 & \sso201011 & 26  \\
(3^3,1^3)     & \sso002000 & \sso200200 & 28  \\
(5,2^2,1^3)   & \sso210100 & \sso111100 & 30 \\
(3^2,2^2,1^2) & \sso010100 & \sso020200 & 40  \\
(3,2^4,1)     & \sso100011 & \sso111200 & 54  \\
(5,1^7)       & \sso220000 & \sso222000 & 56 \\
(3^2,1^6)     & \sso020000 & \sso221011 & 58  \\
(3,2^2,1^5)   & \sso101000 & \sso211111 & 66  \\
(2^6)_1       & \sso000020 & \sso020202 & 70  \\
(2^6)_2       & \sso000002 & \sso020220 & 70  \\
(2^4,1^4)     & \sso000100 & \sso202022 & 76  \\
(3,1^9)       & \sso200000 & \sso222200 & 120 \\
(2^2,1^8)     & \sso010000 & \sso222022 ? & 124 \\
(1^{12})      & \sso000000 & \sso222222 & 220  
\end{array}
$
\vspace{1.3cm}

\end{minipage}
\hspace{10pt}        
\begin{minipage}[t]{8cm}
\centerline {\large\bf so$_{14}$}
\vspace{10pt}
 
$
\begin{array}{c|c|c|c}
\mbox{Partition} & \mbox{Diagram} & \mbox{Weight} &  \\ \hline
(13,1)        & \ssso2222222 & \ssso0000000 & 0  \\
(11,3)        & \ssso2222022 & \ssso0100000 & 2  \\
(11,1^3)      & \ssso2222200 & \ssso2000000 & 4  \\
(9,5)         & \ssso2202022 & \ssso0001000 & 4  \\
(9,3,1^2)     & \ssso2220200 & \ssso1010000 & 6  \\
(7^2)         & \ssso0202022 & \ssso0000011 & 6  \\
(7,5,1^2)     & \ssso2020200 & \ssso1000100 & 8  \\
(7,3^2,1)     & \ssso2200200 & \ssso0020000 & 12 \\
(6^2,1^2)     & \ssso0202011 & \ssso0100011 & 12 \\
(9,2^2,1)     & \ssso2221011 & \ssso1110000 & 14 \\
(5^2,3,1)     & \ssso0200200 & \ssso0010100 & 14 \\
(9,1^5)       & \ssso2222000 & \ssso2200000 & 20 \\
(7,3,2^2)     & \ssso2201011 & \ssso0201000 & 20 \\
(5^2,2^2)     & \ssso0201011 & \ssso0200011 & 22 \\
(7,3,1^4)     & \ssso2202000 & \ssso2101000 & 22 \\
(5,4^2,1)     & \ssso1011011 & \ssso1001100 & 22 \\
(5^2,1^4)     & \ssso0202000 & \ssso2100011 & 24 \\
(5,3^3)       & \ssso2000200 & \ssso0110100 & 24 \\
(5,3^2,1^3)   & \ssso2002000 & \ssso2002000 & 28 \\
(4^2,3^2)     & \ssso0101011 & \ssso0101011 & 28  \\
(7,2^2,1^3)   & \ssso2210100 & \ssso1111000 & 30 \\
(4^2,3,1^3)   & \ssso0110100 & \ssso2000200 & 32  \\
(5,3,2^2,1^2) & \ssso2010100 & \ssso0202000 & 40 
\end{array}
$
\end{minipage}
        
\newpage
\begin{minipage}[t]{8.6cm}
\centerline {\large\bf so$_{14}$}
\vspace{7pt}
$
\begin{array}{c|c|c|c}
\mbox{Partition} & \mbox{Diagram} & \mbox{Weight} & \\ \hline
(3^4,1^2)     & \ssso0002000 & \ssso1011011 & 40  \\
(4^2,2^2,1^2) & \ssso0200011 & \ssso0201011 & 42  \\
(5,2^4,1)     & \ssso2100011 & \ssso1112000 & 54 \\
(7,1^7)       & \ssso2220000 & \ssso2220000 & 56 \\
(3^3,2^2,1)   & \ssso0010100 & \ssso1111011 & 56  \\
(5,3,1^6)     & \ssso2020000 & \ssso2210100 & 58 \\ 
(4^2,1^6)     & \ssso0201000 & \ssso2201011 & 62  \\
(3^3,1^5)     & \ssso0020000 & \ssso2200200 & 64  \\
(5,2^2,1^5)   & \ssso2101000 & \ssso2111100 & 66 \\
(3^2,2^4)     & \ssso0100011 & \ssso0202011 ? & 70 \\
(2^6,1^2)     & \ssso0000011 & \ssso0202022 ? & 112 \\
(5,1^9)       & \ssso2200000 & \ssso2222000 & 120 \\
(3^2,1^8)     & \ssso0200000 & \ssso2221011 & 122 \\
(2^4,1^6)     & \ssso0001000 & \ssso2202022 ? & 140 \\
(3^2,2^2,1^4) & \ssso0101000 & & \\
(3,2^4,1^3)   & \ssso1000100 & & \\
(3,2^2,1^7)   & \ssso1010000 & & \\
(3,1^{11})    & \ssso2000000 & \ssso2222200 & 220 \\
(2^2,1^{10})  & \ssso0100000 & \ssso2222022 ? & 224 \\
( 1^{14})     & \ssso0000000 & \ssso2222222 & 364 
\end{array}
$

\vspace{0.5cm}
\centerline {\large\bf so$_{16}$}
\vspace{7pt}
$
\begin{array}{c|c|c|c}
\mbox{Partition} & \mbox{Diagram} & \mbox{Weight} & \\ \hline
(15,1)        & \ssoo22222222 & \ssoo00000000 & 0  
\end{array}
$
\end{minipage}       
\begin{minipage}[t]{8.6cm}
\centerline{\large\bf so$_{16}$}
\vspace{7pt}
$
\begin{array}{c|c|c|c}
\mbox{Partition} & \mbox{Diagram} & \mbox{Weight} & \\ \hline
(13,3)        & \ssoo22222022 & \ssoo01000000 & 2  \\
(13,1^3)      & \ssoo22222200 & \ssoo20000000 & 4  \\
(11,5)        & \ssoo22202022 & \ssoo00010000 & 4  \\
(11,3,1^2)    & \ssoo22220200 & \ssoo10100000 & 6  \\
(9,7)         & \ssoo20202022 & \ssoo00000100 & 6  \\
(9,5,1^2)     & \ssoo22020200 & \ssoo10001000 & 8  \\
(8^2)_1       & \ssoo02020202 & \ssoo00000002 & 8  \\
(8^2)_2       & \ssoo02020220 & \ssoo00000020 & 8  \\
(7^2,1^2)     & \ssoo02020200 & \ssoo10000011 & 10 \\
(9,3^2,1)     & \ssoo22200200 & \ssoo00200000 & 12 \\
(11,2^2,1)    & \ssoo22221011 & \ssoo11100000 & 14 \\
(7,5,3,1)     & \ssoo20200200 & \ssoo00101000 & 14 \\
(6^2,3,1)     & \ssoo02011011 & \ssoo00010100 & 18  \\
(11,1^5)      & \ssoo22222000 & \ssoo22000000 & 20 \\
(9,3,2^2)     & \ssoo22201011 & \ssoo02010000 & 20 \\
(5^3,1)       & \ssoo00200200 & \ssoo00002000 & 20  \\
(9,3,1^4)     & \ssoo22202000 & \ssoo21010000 & 22 \\
(7,5,2^2)     & \ssoo20201011 & \ssoo02000100 & 22 \\
(7,4^2,1)     & \ssoo21011011 & \ssoo10011000 & 22 \\
(7,5,1^4)     & \ssoo20202000 & \ssoo21000100 & 24 \\
(7,3^3)       & \ssoo22000200 & \ssoo01101000 & 24 \\
(6^2,2^2)_1   & \ssoo02020002 & \ssoo02000002 & 24  \\
(6^2,2^2)_2   & \ssoo02020020 & \ssoo02000020 & 24  
\end{array}
$
\end{minipage}

\newpage
\begin{minipage}[t]{8.73cm}
\centerline {\large\bf so$_{16}$}
\vspace{7pt}
$
\begin{array}{c|c|c|c}
\mbox{Partition} & \mbox{Diagram} & \mbox{Weight} & \\ \hline
(5^2,3^2)     & \ssoo02000200 & \ssoo01100011 & 26  \\
(7,3^2,1^3)   & \ssoo22002000 & \ssoo20020000 & 28 \\
(6^2,1^4)     & \ssoo02020100 & \ssoo20100011 & 28  \\
(5^2,3,1^3)   & \ssoo02002000 & \ssoo20010100 & 30  \\
(9,2^2,1^3)   & \ssoo22210100 & \ssoo11110000 & 30 \\
(5,4^2,3)     & \ssoo10101011 & \ssoo00110011 & 34  \\
(5,4^2,1^3)   & \ssoo10110100 & \ssoo11001100 & 38  \\
(7,3,2^2,1^2) & \ssoo22010100 & \ssoo02020000 & 40 \\ 
(5,3^3,1^2)   & \ssoo20002000 & \ssoo10110100 & 40  \\
(4^4)_1       & \ssoo00020002 & \ssoo00020002 & 40 \\
(4^4)_2       & \ssoo00020020 & \ssoo00020020 & 40 \\
(5^2,2^2,1^2) & \ssoo02010100 & \ssoo02010100 & 42  \\
(7,2^4,1)     & \ssoo22100011 & \ssoo11120000 & 54 \\
(9,1^7)       & \ssoo22220000 & \ssoo22200000 & 56 \\ 
(5,3^2,2^2,1) & \ssoo20010100 & \ssoo11110100 & 56  \\
(7,3,1^6)     & \ssoo22020000 & \ssoo22101000 & 58 \\
(4^2,3,2^2,1) & \ssoo01100011 & \ssoo11101011 & 60 \\
(5^2,1^6)     & \ssoo02020000 & \ssoo22100011 & 60  \\
(3^5,1)       & \ssoo00002000 & \ssoo00200200 & 60 \\
(5,3^2,1^5)   & \ssoo20020000 & \ssoo22002000 & 64 \\
(7,2^2,1^5)   & \ssoo22101000 & \ssoo21111000 & 66 \\
(5,3,2^4)     & \ssoo20100011 & \ssoo02020100 ? & 70 \\
(4^2,2^4)_1   & \ssoo02000002 & \ssoo02020002 & 72 
\end{array}
$
\end{minipage}
\begin{minipage}[t]{8.73cm}
\centerline{\large\bf so$_{16}$}
\vspace{7pt}
$
\begin{array}{c|c|c|c}
\mbox{Partition} & \mbox{Diagram} & \mbox{Weight} & \\ \hline
(4^2,2^4)_2   & \ssoo02000020 & \ssoo02020020 & 72 \\
(3^4,1^4)     & \ssoo00020000 & \ssoo21011011 & 76\\
(4^2,2^2,1^4) & \ssoo02000100 & \ssoo20200200 & 80 \\
(7,1^9)       & \ssoo22200000 & \ssoo22220000 & 120 \\
(5,3,1^8)     & \ssoo20200000 & \ssoo22210100 & 122 \\
(4^2,1^8)     & \ssoo02010000 & \ssoo22201011 & 126 \\
(3^3,1^7)     & \ssoo00200000 & \ssoo22200200 & 128 \\
(2^8)_1       & \ssoo00000002 & \ssoo02020202 ? & 168 \\ 
(2^8)_2       & \ssoo00000020 & \ssoo02020220 ? & 168 \\
(5,1^{11})    & \ssoo22000000 & \ssoo22222000 ? & 220 \\
(3^2,1^{10})  & \ssoo02000000 & \ssoo22221011 & 222 \\
(3,1^{13})    & \ssoo20000000 & \ssoo22222200 & 364 \\
(5,3,2^2,1^4) & \ssoo20101000 & & \\
(5,2^4,1^3)   & \ssoo21000100 & & \\
(5,2^2,1^7)   & \ssoo21010000 & & \\ 
(4^2,3^2,1^2) & \ssoo01010100 & & \\
(4^2,3,1^5)   & \ssoo01101000 & & \\
(3^4,2^2)     & \ssoo00010100 & & \\
(3^3,2^2,1^3) & \ssoo00101000 & & \\
(3^2,2^4,1^2) & \ssoo01000100 & & \\
(3^2,2^2,1^6) & \ssoo01010000 & & \\
(3,2^6,1)     & \ssoo10000011 & & \\
(3,2^4,1^5)   & \ssoo10001000 & & 
\end{array}
$
\end{minipage}
\newpage
        
\begin{minipage}[t]{8.6cm}
\centerline {\large\bf so$_{16}$}
\vspace{7pt}
$
\begin{array}{c|c|c|c}
\mbox{Partition} & \mbox{Diagram} & \mbox{Weight} & \\ \hline
(3,2^2,1^9)   & \ssoo10100000 & & \\
(2^6,1^4)     & \ssoo00000100 & & \\
(2^4,1^8)     & \ssoo00010000 & & \\
(2^2,1^{12})  & \ssoo01000000 & & \\
(1^{16})      & \ssoo00000000 & \ssoo22222222 & 560
\end{array}
$
\end{minipage}
\newpage

\centerline {\bf The Exceptional groups}
\bigskip

\begin{minipage}[t]{5.45 cm}
\vspace{20pt}

\centerline {\bf G$_2$}
\vspace{7pt}
$
\begin{array}{c|c|c}
\mbox{Diagram} & \mbox{Weight} & \\ \hline
\g22 & \g00 & 0  \\
\g02 & \g10 & 2  \\
\g10 & \g11 & 14 \\
\g01 & \g30 & 18 \\
\g00 & \g22 & 56 
\end{array}
$

\vspace{1.5cm}
\centerline {\bf F$_4$}
\vspace{7pt}
$
\begin{array}{c|c|c}
\mbox{Diagram} & \mbox{Weight} & \\ \hline
\f2222 & \f0000 & 0   \\
\f2202 & \f0001 & 2   \\
\f0202 & \f0010 & 6   \\
\f2200 & \f0002 & 8   \\
\f0200 & \f0100 & 12   \\
\f1012 & \f2000 & 16   \\
\f1010 & \f1100 & 28  \\
\f0101 & \f0110 & 34  \\
\f2001 & \f1003 & 34  \\
\f0010 & \f0021 & 38  \\
\f2000 & \f0022 & 56  \\
\f0002 & \f2101 & 70  \\
\f0100 & \f2020 & 72  \\
\f0001 & \f2103 & 118  \\
\f1000 & \f1122 & 156  \\
\f0000 & \f2222 & 312 
\end{array}
$

\end{minipage}
\hspace{5pt}
\begin{minipage}[t]{5.45cm}
\centerline{\bf E$_6$}
\vspace{7pt}
$
\begin{array}{c|c|c}
\mbox{Diagram} & \mbox{Weight} & \\ \hline
\e222222 & \e000000 & 0  \\ 
\e222022 & \e010000 & 2  \\
\e220202 & \e100001 & 4  \\
\e200202 & \e000100 & 6  \\
\e121011 & \e110001 & 10 \\
\e111011 & \e001010 & 12 \\
\e211012 & \e010100 & 14 \\
\e020200 & \e200002 & 16  \\ 
\e000200 & \e100101 & 18  \\ 
\e220002 & \e120001 & 20  \\
\e011010 & \e110101 & 30  \\
\e100101 & \e001110 & 34  \\
\e120001 & \e220002 & 40  
\end{array}
$
\vspace{0.5cm}

\end{minipage}
\hspace{5pt}
\begin{minipage}[t]{6cm}
\centerline{\bf E$_6$}
\vspace{7pt}
$
\begin{array}{c|c|c}
\mbox{Diagram} & \mbox{Weight} & \\ \hline
\e001010 & \e111011 & 42  \\
\e200002 & \e020200 & 56  \\
\e110001 & \e121011 & 60  \\
\e020000 & \e211012 & 70  \\
\e000100 & \e111111 & 78  \\
\e100001 & \e220202 & 120 \\
\e010000 & \e222022 & 168 \\
\e000000 & \e222222 & 312 \\  
\end{array}
$

\vspace{0.7cm}
\centerline {\bf E$_7$}
\vspace{7pt}
$
\begin{array}{c|c|c}
\mbox{Diagram} & \mbox{Weight} & \\ \hline
\ee2222222 & \ee0000000 & 0  \\
\ee2220222 & \ee1000000 & 2  \\
\ee2220202 & \ee0000010 & 4  \\
\ee2002022 & \ee0010000 & 6  
\end{array}
$
\end{minipage}

\newpage

\hspace{50pt}
\begin{minipage}[t]{6cm}
\centerline {\bf E$_7$}
\vspace{7pt}
$
\begin{array}{c|c|c}
\mbox{Diagram} & \mbox{Weight} & \\ \hline
\ee2022020 & \ee0000002 & 6  \\
\ee2002020 & \ee0100001 & 8  \\
\ee2002002 & \ee0001000 & 12 \\
\ee2110122 & \ee1010000 & 14 \\
\ee0002020 & \ee0200000 & 14 \\
\ee2110110 & \ee0000020 & 16 \\
\ee0002002 & \ee0010010 & 18 \\ 
\ee2110102 & \ee2000010 & 20 \\
\ee2020020 & \ee2000002 & 22 \\
\ee0020020 & \ee1000101 & 24 \\
\ee2000200 & \ee0001010 & 28 \\
\ee0110102 & \ee1010010 & 30 \\
\ee0002000 & \ee0000200 & 30 \\ 
\ee1001020 & \ee1100100 & 32 
\end{array}
\bigskip
$
\end{minipage}
\hspace{50pt}
\begin{minipage}[t]{6cm}
\centerline{\bf E$_7$}
\vspace{7pt}
$
\begin{array}{c|c|c}
\mbox{Diagram} & \mbox{Weight} & \\ \hline
\ee1001012 & \ee0011000 & 34 \\
\ee2001010 & \ee2000020 & 40 \\
\ee1001010 & \ee1001010 & 42 \\
\ee0000200 & \ee0002000 & 48 \\ 
\ee2110001 & \ee1010020 & 54 \\
\ee2000022 & \ee2020000 & 56 \\
\ee2000020 & \ee2110001 & 58  \\
\ee0001010 & \ee2000200 & 62  \\
\ee0110001 & \ee1001020 & 70  \\
\ee2020000 & \ee2000022 & 70  \\ 
\ee0020000 & \ee1001012 & 72  \\
\ee1000101 & \ee1011010 & 78  \\
\ee1001000 & \ee1101011 & 84  \\
\ee0010010 & \ee0011101 & 88  
\end{array}
$
\end{minipage}

\newpage

\hspace{40pt}
\begin{minipage}[t]{7cm}
\vspace{10pt}
\centerline{\bf E$_7$}
\vspace{7pt}
$
\begin{array}{c|c|c}
\mbox{Diagram} & \mbox{Weight} & \\ \hline
\ee0200000 & \ee0002020 & 112 \\
\ee2000002 & \ee2020020 & 120 \\
\ee0000020 & \ee2110110 & 122 \\
\ee2000010 & \ee2110102 & 124 \\
\ee0001000 & \ee2002002 & 126 \\
\ee1000010 & \ee2110120 & 166 \\
\ee0100001 & \ee2002020 ? & 168 \\
\ee2000000 & \ee2110122 & 220 \\  
\ee0010000 & \ee2002022 ? & 222 \\
\ee0000002 & \ee2022020 & 312  \\
\ee0000010 & \ee2220202 ? & 318 \\
\ee1000000 & \ee2220222 ? & 462 \\
\ee0000000 & \ee2222222 & 798 
\end{array}
$

\end{minipage}
\hspace{20pt}
\begin{minipage}[t]{7cm}
\centerline{\bf E$_8$}
\vspace{7pt}
$
\begin{array}{c|c|c}
\mbox{Diagram} & \mbox{Weight} & \\ \hline
\eee22222222 & \eee00000000 & 0 \\
\eee22202222 & \eee00000001 & 2  \\
\eee22202022 & \eee10000000 & 4  \\
\eee20020222 & \eee00000010 & 6  \\
\eee20020202 & \eee01000000 & 8  \\
\eee20020022 & \eee00000100 & 12 \\
\eee21101222 & \eee00000011 & 14 \\
\eee20020020 & \eee00100000 & 14 \\
\eee21101022 & \eee10000002 & 20 \\
\eee21101101 & \eee11000000 & 22 \\
\eee20002002 & \eee10000100 & 28 \\
\eee01101022 & \eee10000011 & 30 \\
\eee20010102 & \eee20000002 & 40 \\
\eee10010102 & \eee10000101 & 42  

\end{array}
$
\end{minipage}

\newpage

\hspace{40pt}
\begin{minipage}[t]{7cm}
\centerline {\bf E$_8$}
\vspace{7pt}
$
\begin{array}{c|c|c}
\mbox{Diagram} & \mbox{Weight} &  \\ \hline
\eee00002002 & \eee00000200 & 48 \\ 
\eee00020022 & \eee00100100 & 50 \\
\eee21100012 & \eee20000011 & 54 \\
\eee10010110 & & \\ 
\eee10010101 & & \\
\eee20000222 & \eee00000022 & 56 \\
\eee20000202 & \eee01000012 & 58 \\
\eee00002000 & \eee00010010 & 60 \\
\eee20000200 & \eee02000002 & 64 \\
\eee10001012 & \eee10000111 & 78 \\
\eee00010100 & & \\ 
\eee10001010 & & \\
\eee01100010 & & \\
\eee00100101 & & 
\end{array}
$
\end{minipage}
\hspace{20pt}
\begin{minipage}[t]{7cm}
\centerline{\bf E$_8$}
\vspace{7pt}
$
\begin{array}{c|c|c}
\mbox{Diagram} & \mbox{Weight} & \\ \hline
\eee10010001 & & \\
\eee00010010 & \eee00002000 ? & 80 \\
\eee00010102 & \eee11001000 & 82 \\ 
\eee00020020 & \eee01100011 & 98 \\
\eee10010100 & & \\
\eee01100012 & \eee10010011 & 116 \\
\eee00020002 & \eee00020010 & 174 \\
\eee02000002 & \eee20000200 & 112 \\ 
\eee20000022 & \eee20000022 & 120 \\
\eee20000020 & \eee10001012 & 122 \\
\eee00100100 & & \\
\eee00010002 & \eee00010102 & 126 \\
\eee20000101 & & \\
\eee00000200 & \eee00002002 & 128  
\end{array}
$
\end{minipage}

\newpage

\hspace{40pt}
\begin{minipage}[t]{7cm}
\centerline {\bf E$_8$}
\vspace{7pt}
$
\begin{array}{c|c|c}
\mbox{Diagram} & \mbox{Weight} & \\ \hline
\eee00010001 & & \\
\eee10001000 & & \\
\eee10000101 & & \\
\eee02000000 & \eee10010110 & 168 \\ 
\eee01000012 & \eee20000202 ? & 168\\
\eee00010000 & & \\
\eee20000002 & \eee21100012 & 220 \\
\eee10010122 & \eee12010002 & 230 \\  
\eee10000100 & & \\
\eee01000010 & & \\
\eee00100001 & & \\ 
\eee00001000 & & \\
\eee00000022 & \eee20000222 & 312 \\
\eee00000020 & \eee10010122 & 314 
\end{array}
$
\end{minipage}
\hspace{25pt}
\begin{minipage}[t]{7cm}
\centerline{\bf E$_8$}
\vspace{7pt}
$
\begin{array}{c|c|c}
\mbox{Diagram} & \mbox{Weight} & \\ \hline
\eee10000102 & \eee21100112 & 330 \\
\eee00000101 & & \\
\eee10000010 & & \\
\eee20000000 & \eee11101110 & 336 \\
\eee00100000 & \eee20020020 ? & 368 \\
\eee10000002 & \eee21101022 ? & 462 \\
\eee00000100 & \eee20020022 ? & 464 \\
\eee10000001 & & \\
\eee01000000 & \eee20020202 ? & 560 \\
\eee00000002 & \eee21101222 & 798  \\
\eee00000010 & \eee20020222 ? & 800 \\
\eee10000000 & \eee22202022 ? & 1040 \\ 
\eee00000001 & \eee22202222 ? & 1520 \\
\eee00000000 & \eee22222222 & 2480 
\end{array}
$
\end{minipage}

\end{document}